\documentclass[12pt]{article}
\usepackage{amssymb}
\usepackage{a4wide}
\usepackage{epsfig}
\usepackage{oldgerm}
\usepackage{amsmath}
\usepackage{cite}
\usepackage{amsthm}
\usepackage{fixmath}
\usepackage{amstext,color}
\usepackage{array,multirow}
\usepackage{fullpage}
\usepackage{graphicx,graphics}

\usepackage[vlined,ruled,linesnumbered]{algorithm2e}

\definecolor{teal}{rgb}{0.0,0.5,0.5}
\definecolor{frenchrose}{rgb}{0.96,0.29,0.54}
\definecolor{lg}{rgb}{0.36,0.99,0.82}
\definecolor{dblue}{rgb}{0,0,.5}
\definecolor{dpink}{cmyk}{.2,1,.1,.04}
\definecolor{purple}{rgb}{0.35,0.04,0.64}
\definecolor{borange}{rgb}{1, .388, 0}
\definecolor{dpurple}{rgb}{0.61,0.22,1.00}
\definecolor{purp}{rgb}{0.44,0.00,0.87}
\definecolor{green}{rgb}{0.00,0.44,0.00}
\definecolor{junebud}{rgb}{0.74,0.85,0.34}
\definecolor{plum}{rgb}{0.56,0.27,0.52}
\newtheorem{theorem}{Theorem}[section]
\newtheorem{corollary}{Corollary}[theorem]

\newtheorem{definition}{Definition}[section]
\newtheorem{lemma}{Lemma}[section]
\newtheorem{remark}{Remark}[section]

\newtheorem{prop}{Proposition}
\def\CC{{\textmd \kern.24em \vrule width.02em height1.4ex depth-.05ex \kern-.26emC}}
\def\TagOnRight
\def\QQ{\rlap {\raise 0.4ex \hbox{$\scriptscriptstyle |$}} {\hskip -0.1em Q}}
\catcode`\@=11

\begin{document}
%%%%%%%%%%%%%%%%%%%%%%%%%%%%%%%%
\begin{center}
 {{\bf \large {\rm {\bf  L1 scheme on graded mesh for subdiffusion equation with nonlocal diffusion term }}}}
\end{center}
\begin{center}
{\textmd {{\bf Sudhakar Chaudhary,}}}\footnote{\it Department of Mathematics,  Institute of Infrastructure, Technology, Research and Management, Ahmedabad, Gujarat, India, (sudhakarchaudhary@iitram.ac.in)}
{\textmd {\bf Pari J. Kundaliya}}\footnote{\it Department of Mathematics,  Institute of Infrastructure, Technology, Research and Management, Ahmedabad, Gujarat, India, (pariben.kundaliya.pm@iitram.ac.in)}
%%\\
%%%{\textmd {{\bf V.V.K Srinivas Kumar,}}}\footnote{\it Department of Mathematics, Indian Institute of Technology, Delhi, India, {(vvksrini@maths.iitd.ac.in)} }
%%%{\textmd {{\bf Balaji Srinivasan.}
%%%}\footnote{{\it Department of Applied Mechanics, Indian Institute of Technology, Delhi, India, } (balaji@am.iitd.ac.in)}}
\end{center}
%%%%%%%%%%%%%%%%%%%%%%%%%%%%%%%%%%%%%%page 1 %%%%%%%%%%%%%%%%%%%%%%%%%%%%%%%%%%%%%%%%%%%%%%%%%
\begin{abstract} The solution of time fractional partial differential equations in general exhibit a weak singularity near the initial time. In this article we propose a method for solving time fractional diffusion equation with nonlocal diffusion term. The proposed method comprises L1 scheme on graded mesh, finite element method and Newton's method. We discuss the well-posedness of the weak formulation at discrete level and derive \emph{a priori} error estimates for fully-discrete formulation in $L^2(\Omega)$ and $H^1(\Omega)$ norms. Finally, some numerical experiments are conducted to validate the theoretical findings.
\end{abstract}
%\begin{center}
{\bf Keywords:}  Nonlocal problem; initial singularity; L1 scheme; graded mesh; error estimate.\\
\noindent {\bf AMS(MOS):} 65M12, 65M60, 35R11.
%\end{center}
%%%%%%%%%%%%%%%%%%%%%%%%%%%%%%%%%%
%
\section{Introduction}
The study of nonlocal problems has gained considerable attentions in recent years (see, \cite{[me],[sk],[ch],[r6new],[rbl]} and references therein).
%%%%%%%%%%%%
Author in \cite{[r6new]} consider the following parabolic nonlocal problem
	\begin{equation}\label{bjn}
	\begin{split}
	u_t- a\big(l(u)\big) \, \Delta u + f(u) &= h  \quad \mbox{in}  \quad  Q = \Omega\times(0,T),\\
	u(x,t)&=0 \quad \mbox{on} \quad \Sigma = \Gamma \times (0,T),\\
	u(x,0)&=u_{0}(x) \quad \mbox{in} \quad \Omega,
	\end{split}
	\end{equation}
	where $\Omega$ is a smooth bounded open subset of $\mathbb{R}^d$ $(d \ge 1)$ with smooth boundary and
$l \, : \, $ $L^2(\Omega) \longrightarrow  \mathbb{R}$ is a continuous linear form.
	This type of problems, besides its mathematical motivation, arises from physical situations related to migration of a population of bacterias in a container (say domain $\Omega$) where $u$ describes the density of the population of bacteria. The velocity of migration ($\vec{v}$) of this population is proportional to the gradient of the density with a {positive} factor $'a'$ depending on entire population, that is $\vec{v} = a \, \nabla u$, $a = a \big( \int_{\Omega} \, u \, dx\big)$.  Authors in \cite{[r4new]} have mentioned that fractional order differential equations are naturally related to systems with memory which exists in most biological systems. Also, in the nature there exist many processes that can not be adequately described with classical exponential law which is corresponding to integer order derivative. For example, the dynamics of population densities can follow a law that behaves like the exponential one but changes slowly or faster than the exponential function (for more details see \cite{[r5new]}).
%%%%%%%%%%%%%%%%%%%%%
%Parabolic partial differential equations (PDEs) are popular in the modeling of anomalous slow transport process and very often these equations are called fractional diffusion or subdiffusion equations. have been applied in various field of science and engineering, such as physics, chemistry, biology and finance \cite{[szy]}. The optimal $L^2- norm$ error estimate and optimal $H^1- norm$ error estimate are proved in \cite{[r02]} and \cite{[hr12]}, respectively for a time-fractional diffusion equation.
 This motivates us to consider following time fractional partial differential equation (PDE) with nonlocal diffusion term: Find $u$ such that
  \begin{subequations}\label{e1}
  	\begin{align}
  	\label{cuc3:1.1}
  	^{c}_{0}{D}^{\alpha}_{t}u(x,t)-a(l(u)) \: \Delta u(x,t)&=f(u)  \quad \mbox{in}  \quad \Omega\times(0,T], \\
  	\label{cuc3:1.2}
   	u(x,t)&=0  \quad \mbox{on}  \quad \partial\Omega\times(0,T],\\
   	\label{cuc3:1.3}
  	u(x,0)&=u_0(x) \quad \mbox{in} \quad \Omega,
  	\end{align}
  \end{subequations}
  where $\Omega$ is an open bounded subset of $\mathbb{R}^d$ $(d=1$ or $2)$ with smooth boundary $\partial \Omega$,  $l(u)= \int_{\Omega} u(x,t) \, dx$ and the Caputo fractional derivative $^{c}_{0}{D}^{\alpha}_{t}u(x,t)$ $(0 < \alpha < 1)$ is defined as \cite{[r1]} \\
  \begin{equation}\label{e2}
  \begin{split}
  {^{c}_{0}D^{\alpha}_{t}}u(x,t) =& \frac{1}{\Gamma (1-\alpha)}\int_{0}^{t}(t-s)^{-\alpha} \; \frac{\partial u(x,s)}{\partial s} \, ds , \;\; t>0.
  \end{split}
  \end{equation}
Problem \eqref{cuc3:1.1}- \eqref{cuc3:1.3} is a time fractional version of the integer order nonlocal parabolic problem given in \cite{[me],[sk],[r6new]}. Authors in  \cite{[mn1],[mn2],[mn3]} have studied similar nonlocal problem using finite element method with uniform time grid.
Due to wide range of applications, the development of effective numerical methods for time-fractional parabolic PDEs is of great importance. Many efficient time-stepping schemes have been proposed in the literature for linear time fractional PDEs. Mainly, these schemes can be divided into two class: L1 type schemes and convolution quadrature (CQ) (See \cite{[r8],[r14],[xzw1],[xzw3],[xzw4],[hr12],[ach1],[r02],[r11],[r5a],[r13],[r10],[r9]} and references therein). An important consideration in the development of numerical methods for fractional diffusion problems is that solution $u$ is weakly singular in time near time $t=0$ \cite{[kwn],[r02],[r08],[r8],[r13],[r9]}. Authors in \cite{[r8]} consider the following time fractional diffusion equation:
\begin{equation}\label{bjn}
	\begin{split}
	^c_{0}{D}^{\alpha}_{t}u(x,t)-\Delta{u(x,t)} &= f(u)  \quad \mbox{in}  \quad \Omega\times(0,T],\\
	u(x,t)&=0 \quad \mbox{on}\quad\partial{\Omega}\times(0,T],\\
	u(x,0)&=u_{0}(x),\quad \mbox{in} \quad \Omega,
	\end{split}
\end{equation}
They show that if $u_0\in H^1_0(\Omega) \cap H^2(\Omega)$ then the solution to problem \eqref{bjn} satisfies $\|\partial_t u\|_{L^2{(\Omega)}}\leq C t^{\alpha-1}.$ The spaces $L^2{(\Omega)}$ and $H^m(\Omega),$ for $m=1, 2$ are defined in the next section.
%In our work also we assume that solution to problem \eqref{e1} is weakly singular $t=0.$
In our work also, we assume that the solution to problem \eqref{e1} satisfies $\|u^{(l)}(t)\|_{L^2{(\Omega)}} \lesssim 1 + t^{\alpha- l}$, for $l=0,1,2$ and $t \in (0,T]$.
In case of weak singularity  near $t=0,$ the $L1$ scheme on uniform mesh for problem \eqref{e1} gives $\mathcal{O}(\tau^\alpha)$ order of convergence in maximum norm in time, where $\tau$ denotes the time step size (see the numerical section) \cite{[r8]}. To overcome this issue, we analyse the L1 scheme on graded mesh \cite{[kwn],[r02]} for the problem \eqref{e1} having initial singularity at $t=0.$ We have shown that with appropriate choice of the grading parameter $r,$ one can recover the optimal convergence order $2-\alpha$ in temporal direction. Another concern of this work is to derive the optimal order of convergence in $L^2(\Omega)$ and $H^1(\Omega)$ norms in spatial direction. Authors in \cite{[r5a]} have pointed that due to initial singularity of the solution and the discrete convolution form in numerical Caputo derivative, the traditional $H^1(\Omega)$- norm analysis (corresponding to the case for a classical diffusion equation) to the time fractional diffusion problem always leads to suboptimal estimates. A similar conclusion is also drawn in \cite{[ddgm]}, where authors have used direct discontinuous Galerkin method for solving the time fractional diffusion equation. For the derivation of optimal error estimate in $H^1(\Omega)$ norm, we follow the idea given in \cite{[r5a],[ach1]}. To the best of our knowledge this is the first attempt when $L1$ scheme on graded mesh is used for solving the subdiffusion equation with nonlocal diffusion term. The main contribution of the present work are summarized below.\\
$\bullet$ To handle the weak singularity in the solution, we approximate $^{c}_{0}{D}^{\alpha}_{t}u(x,t)$ by well known $L1$ scheme on graded temporal mesh \cite{[r02]}.\\
$\bullet$ For nonlocal term and nonlinearity in right-hand side, we use Newton's method.\\
$\bullet$ We derive a priori bound for the fully discrete solution in $L^2(\Omega),$ $H_0^1(\Omega)$ norms and prove the existence-uniqueness of fully discrete solution.\\
$\bullet$ In previous papers \cite{[mn1],[mn2],[mn3]}, the authors have shown convergence in $L^2(\Omega)$ norm using uniform mesh. We prove optimal rate of convergence in $L^2(\Omega)$ and $H^1(\Omega)$ norms using graded mesh. \\Throughout the paper, $C$ denote a generic constant (not necessarily same at different occurrences) while $C_i,$ $i=1,2,3 ... 15,$ are fixed constants; all these constants are independent of mesh parameters $h$ and $N.$

The rest of the paper is organized as follows: In Section 2, first we recall some basic definitions and then define the weak formulation of given problem \eqref{e1}. In Section 3, we give fully-discrete scheme and derive the \emph{a priori} bound for the fully-discrete scheme. We prove the existence-uniqueness of fully-discrete solution in section 4. Error analysis of our proposed scheme is presented in Section 5. Finally, numerical results in Section 6 confirm our theoretical estimates.
%%%%%%%%%%%%%%%%%%%%%%%%%%%%%%%%%%

\section{Preliminaries and weak formulation}
 Let $L^2(\Omega)$ be the space of square integrable functions on $ \Omega $ with inner product $(g_1,g_2) = \int_{\Omega} \, g_1(x) g_2(x) \, dx$ and norm  $ \|{g_1}\| = {\Big(\int_{\Omega} \, |g_1(x)|^2 \, dx\Big)}^{\frac{1}{2}}.$ For a non-negative integer $m,$ $H^m(\Omega)$ denote the Sobolev space on $\Omega$ with the norm $\|w\|_m=\Big(\sum_{0\leq \alpha \leq m}\|\frac{\partial^{\alpha}w}{\partial x^{\alpha}}\|^2\Big)^{1/2}.$ Also, we write $\|w\|_{L^{\infty}(H^m)} = \sup_{0<t\leq T}\|w\|_m$ and \[ H^1_0(\Omega) := \Big\{ v \in H^1(\Omega) \, : v \, = \, 0 \: \mbox{on} \: \partial \Omega   \Big\} . \]\\
%\noindent Again for $w:[0,T]\rightarrow H^m(\Omega),$ define the norms $\|w\|_{L^2(0,T;H^m(\Omega))},$ $||w||_{L^{\infty}(0,T;L^2(\Omega))}$ and $||w||_{L^{\infty}(0,T;H_0^1(\Omega))}$ by
%\begin{equation*}
%\begin{split}
% \|w\|_{L^2(0,T;H^m(\Omega))}=&\left(\int_{0}^{T}\|w(t)\|_m^2\right)^{1/2},\\
%\|w\|_{L^{\infty}(0,T;L^2(\Omega))}=&ess\sup_{0\leq t\leq T}\|w(t)\|,\\
%\|w\|_{L^{\infty}(0,T;H_0^1(\Omega))}=&ess\sup_{0\leq t\leq T}\|\nabla w(t)\|.
%\end{split}
%\end{equation*}
 The {weak formulation} of given problem \eqref{e1} is: find $u(\cdot,t) \in  H^1_0(\Omega)$ such that for each $t \in(0,T]$ one has
 \begin{equation}\label{e3}
 \begin{split}
   (^{c}_{0}{D}^{\alpha}_{t}u, v ) \, + \, a\big(l(u)\big) \, (\nabla u, \nabla v) \, =& \, \big(f(u), v\big), \quad  \forall v \in H^1_0(\Omega).\\
 u(x,0)=& \, u_0(x), \quad \mbox{in} \; \, \Omega.
 \end{split}
 \end{equation}
 In our further analysis, we need following hypotheses on given data.
 \begin{itemize}
 	\item {H1:} $a: \mathbb{R} \rightarrow \mathbb{R}$ is bounded with  $0<m_1 \le a(x) \le m_2 < \infty, \; \forall \, x \in \mathbb{R}.$\\
 	  \item {H2:} $a: \mathbb{R} \rightarrow \mathbb{R}$ is Lipschitz continuous with Lipschitz constant $L>0,$ \textit{i.e.,}
 	\begin{equation}\label{e5}
 	|a(x_1)-a(x_2)| \le L\, |x_1-x_2|, \quad \forall \, x_1, x_2 \in \mathbb{R}.
 	\end{equation}
 	%\item {H3:} $l: L^2(\Omega) \rightarrow \mathbb{R}$ is continuous linear form, \textit{i.e} $\exists \, g \in L^2(\Omega)$ such that \\ \indent \ $l(u) = \int_{\Omega} \, g(x) \, u(x) \:dx$ and
% 	\begin{equation}\label{e6}
% 	|l(u)| \le C \, \|{u}\|_{L^2(\Omega)}, \; \forall \, u \in L^2(\Omega).
% 	\end{equation}
 	 \item {H3:} $u_0 \in H_0^1(\Omega)\cap H^2(\Omega)$ and $f : \mathbb{R} \rightarrow \mathbb{R}$ is Lipschitz continuous with Lipschitz constant $K>0,$ \textit{i.e.,} $ |f(u_1) - f(u_2)| \leq K|u_1 - u_2|$, for $u_1,u_2$ $\in$ $\mathbb{R}$.

 \end{itemize}
	Under the hypotheses H1-H3, it can be shown that there exist a unique weak solution of the problem \eqref{e3}. For the proof, one can use Faedo-Galerkin method in a standard way (\cite{[rbl],[mn1],[mn3],[r08]}). Note that for the derivation of error estimates, we require some additional regularity of solution $u$ and it is mentioned in Section 5.
%%%%%%%%%%%%%%%%%%%%%%%%%%%%%%%%%%

\section{Fully-discrete formulation and \emph{a priori} bound}
For the spatial discretization, we use Galerkin finite element method (FEM). Let $\Omega_h$ be a quasi uniform partition of $\Omega$ into disjoint subintervals in 1D and triangles in 2D with  step size $h$. Let $M>0$ be the finite integer and $X_h$ be the $M$-dimensional subspace of $H^{1}_{0}(\Omega)$ such that $X_h$ consists of continuous functions on closure $\overline{\Omega}$ of $\Omega$ which are linear in each $T_k \, \in \, \Omega_h$ and vanishes on $\partial \Omega$. \textit{i.e.,}
 \begin{equation*}
 X_h:=\Big\{v\in C^{0}(\bar{\Omega}): v_{|{\small T_{k}}}\in P_1(T_k), \: \forall \: T_{k}\in \Omega_h \: \, \mbox{and} \: \, v=0 \: \mbox{on} \: \partial \Omega\Big\}.
 \end{equation*}
 For each $i= 1,2,...,M$, let $\phi _i(x)$ be the pyramid function in $X_h$ which takes the value 1 at $i-th$ node and vanishes at other node points. Then $\left\lbrace \phi _i  \right\rbrace ^{M}_{i=1}$ forms a basis for the space $X_h.$ \\
  For temporal discretization, we use L1 scheme on graded mesh \cite{[kwn],[r02]}. For this let $N \in \mathbb{N}$ and $0=t_0<t_1<t_2 <...<t_N=T$ be a partition of the time interval $[0,T]$ such that $t_n=T(\frac{n}{N})^{r}$, for $n=0,1,...,N$, where $r \ge 1$ is the mesh grading parameter (for $r=1$, the mesh is uniform). Also, the step size is given by $\tau_n = t_n - t_{n-1}$, for $n=1,2,...,N$. For $1 \le n \le N$, let $U^n$ denote the approximate value of $u$ at the node point $t_n$.\\
  Now, the L1-approximation to Caputo fractional derivative on the graded mesh \cite{[kwn],[r02],[r11]} is given below.
\begin{equation}\label{e8}
\begin{split}
^c_0D^{\alpha}_{ t_n}u \, =& \, \frac{1}{\Gamma {(1- \alpha)}} \, \int_{0}^{t_n} (t_n -s)^{- \alpha} \, \frac{ \partial u(x,s)}{ \partial s} \: ds\\
=& \, \frac{1}{\Gamma {(1- \alpha)}} \, \sum_{k=0}^{n-1} w_{n,k} \, \left( u(x,t_{k+1}) - u(x,t_{k}) \right)  + E_n, \quad \forall n = 1,2,...,N, \\
\end{split}
\end{equation}
where
\begin{equation}\label{e9}
\begin{split}
\quad w_{n,k} = \frac{1}{ \tau _{k+1} } \, \int_{t_{k-1}}^{t_k} (t_n -s)^{- \alpha} \, ds, \quad \forall k = 1,2,...,n.
\end{split}
\end{equation}
For any sequence $\left\lbrace v_n \right\rbrace^{N}_{n=1} $, let us define a discrete Caputo fractional differential operator $D^{\alpha} _{N}$ as
\begin{equation}\label{e10}
\begin{split}
D^{\alpha} _{N} v^n \, :=& \, \frac{1}{\Gamma {(2- \alpha)}} \, \sum_{k=0}^{n-1} w_{n,k} \, \left( v^{k+1} - v^{k} \right), \quad \forall n = 1,2,...,N, \\
=& \,\frac{ d_{n,1}}{\Gamma {(2- \alpha)}} \, v^n - \frac{d_{n,n}}{\Gamma {(2- \alpha)}} \, v^0 + \frac{1}{\Gamma {(2- \alpha)}} \, \sum_{k=1}^{n-1} (d_{n,k+1} - d_{n,k}) \, v^{n-k},\\
\end{split}
\end{equation}
where
\begin{equation}\label{e11}
\begin{split}
 d_{n,k} \, = \, \frac{(t_n - t_{n-k})^{1- \alpha} - (t_n - t_{n-k+1})^{1- \alpha}}{\tau_{n-k+1}}, \quad \mbox{for} \; \, 1 \le k \le n \le N.
\end{split}
\end{equation}
In particular, $d_{n,1} \, = \, \tau_n^{-\alpha}$. From mean value theorem, it follows that
\begin{equation}\label{e100}
	d_{n, k+1} \, \le \, d_{n, k}, \quad \mbox{for} \; 0 \le k \le n-1 \le N-1.
\end{equation}
Let $\zeta ^n \, := \, ^c_0D^{\alpha}_{t_n}u(x, t_n) - D^{\alpha} _{N} u(x, t_n)$.
\begin{lemma}\label{l6}
  \cite{[r02], [hr12], [ach1]}  Assume that $|u^{(q)}(t)| \lesssim 1 + t^{\alpha - q}$, for $q=0,1,2$ and $t \in (0,T]$. Then there exists a constant $C$ such that for $n=1, 2, ..., N$,
  %\begin{equation}\label{e99}
%  \| \zeta ^n \| \, \le \, C \, n^{-min \left\lbrace 2- \alpha, \, r \alpha \right\rbrace }, \quad
%  \end{equation}
  \begin{equation}\label{ee99}
  \| \zeta ^n \|_{1} \, \le \, C \, n^{-min \left\lbrace 2- \alpha, \, r \alpha \right\rbrace }.
  \end{equation}	
\end{lemma}

 Now, following the idea given in \cite{[ach1]}, we define coefficients $b^{(n)}_{n-i}$ and $p^{(n)}_{n-i}$ as follows: for $1 \le i \le n \le N$,
 \begin{equation}\label{e101}
 \quad b^{(n)}_{n-i} := \frac{d_{n, \, n-i+1}}{\Gamma{(2- \alpha )}}, \; \; \mbox{for} \; \, 1 \le i \le n \le N.
 \end{equation}
and
 \begin{equation}\label{e102}
   \begin{split}
      	p^{(n)}_{n-i} :=  \left\{\begin{array}{l} \Gamma {(2-\alpha)} \, \tau_i^{\alpha} \sum_{k=i+1}^{n} \, \left( b^{(k)}_{k-i-1} - b^{(k)}_{k-i} \right) \, p^{(n)}_{n-k}, \quad \mbox{if} \; \, i=1, 2,...,(n-1), \\
      	\Gamma {(2-\alpha)} \, \tau_n^{\alpha}, \quad \mbox{if} \; \, i=n. \\
       \end{array}\right.
   \end{split}
 \end{equation}
\begin{lemma}\label{l3}
   \cite{[ach1], [ddgm]} For $n=1,2,...,N$, one has
   \begin{equation}\label{e104}
    \begin{split}
      \sum_{s=1}^{n} \, p^{(n)}_{n-s} \, s^{-min\left\lbrace 2-\alpha, \, r \alpha \right\rbrace } \, \le \, C \, N^{-min\left\lbrace 2-\alpha, \, r \alpha \right\rbrace }.
   \end{split}
   \end{equation}	
\end{lemma}

\begin{lemma}\label{l7}
	\cite{[r5]} For $n=1,2,...,N$, one has
	\begin{equation}\label{e124}
	\begin{split}
	\sum_{s=1}^{n} \, p^{(n)}_{n-s} \, \frac{t^{- \alpha}_s}{\Gamma {(1- \alpha)}} \, \le \, 1.\\
	\end{split}
	\end{equation}	
\end{lemma}

Using the above notations, the fully-discrete scheme for the given problem is as follows: For each $1 \le n \le N$, find $U^n \in X_h$ such that
\begin{equation}\label{e13}
\begin{split}
\left( {D}^{\alpha}_{N}U^n, v_h \right) \, + \, a\big(l(U^n)\big) \, (\nabla U^n, \nabla v_h) \, =& \, \big(f(U^n), v_h\big), \quad  \forall v_h \in X_h, \\
U^0 \, =& \, u_h^0,
\end{split}
\end{equation}
where $u_h^0$ is some approximation of $u_0(x)$.\\
 From \eqref{e13}, we have for $1 \le i \le M,$
\begin{equation}\label{e13n}
\begin{split}
\left( {D}^{\alpha}_{N}U^n, \, \phi_i \right) \, + \, a\big(l(U^n)\big) \, (\nabla U^n, \, \nabla \phi_i) \, = \, \big(f(U^n), \, \phi_i\big).\\
\end{split}
\end{equation}
From \eqref{e10} and \eqref{e13n}, we get
\begin{equation}\label{e47}
\begin{split}
( U^n, \, \phi_i) + \frac{1}{d_{n,1}} \left(- d_{n,n} U^0 + \sum_{k=1}^{n-1} (d_{n,k+1} - d_{n,k}) U^{n-k}, \, \phi_i \right) + \frac{\Gamma{(2- \alpha)}}{d_{n,1}} \, a\big(l(U^n)\big)\\
(\nabla U^n, \, \nabla \phi_i) \,= \, \frac{\Gamma{(2- \alpha)}}{d_{n,1}} \big(f(U^n), \, \phi_i\big),\\
\end{split}
\end{equation}
where weights $d_{n,k+1} \, (1 \le k \le n \le N)$ are given in \eqref{e11}.\\
Since $U^n \in X_h, \,  \exists \ \alpha^n_j \in \mathbb{R}$ such that
\begin{equation}\label{e48}
U^n \,= \, \sum_{j=1}^{M}  \alpha^n_j \phi_j.
\end{equation}
Set $\bar{\alpha}^n \, = [ \alpha^n_1, \alpha^n_2,...,\, \alpha^n_M].$ Putting above value of $U^n$ in \eqref{e47}, we get system of non-linear algebraic equations
\begin{equation}\label{e49}
F_i(\bar{\alpha}^n) \, = \, F_i(U^n) \, = \, 0, \quad \mbox{for} \; 1 \le i \le M,
\end{equation}
where
\begin{equation}\label{e50}
\begin{split}
F_i(U^n) \, = \, &( U^n, \, \phi_i) + \frac{1}{d_{n,1}} \left(- d_{n,n} U^0 + \sum_{k=1}^{n-1} (d_{n,k+1} - d_{n,k}) U^{n-k}, \, \phi_i \right)\\
&+ \frac{\Gamma{(2- \alpha)}}{d_{n,1}} \, a\big(l(U^n)\big) \, (\nabla U^n, \, \nabla \phi_i) \, - \, \frac{\Gamma{(2- \alpha)}}{d_{n,1}} \, \big(f(U^n), \, \phi_i\big).
\end{split}
\end{equation} 	
If we use Newton's method for finding $ \alpha^n_j, \, 1 \le j \le M $ then elements of the Jacobian matrix $J_1$ take the following form
\begin{equation}\label{e51}
\begin{split}
\frac{\partial F_i}{\partial \alpha^n_j} (U^n)  \, = \, &( \phi_j, \, \phi_i) - \, \frac{\Gamma{(2- \alpha)}}{d_{n,1}} \, \big(f'(U^n) \phi_j, \, \phi_i\big) + \frac{\Gamma{(2- \alpha)}}{d_{n,1}} \, a\big(l(U^n)\big) \, (\nabla \phi_j, \, \nabla \phi_i) \\
&+ \frac{\Gamma{(2- \alpha)}}{d_{n,1}} \, a'\big(l(U^n)\big) \: l'(U^n) \, \left( \int_{\Omega} \phi_j \, dx \right) \, (\nabla U^n, \, \nabla \phi_i).\\
\end{split}
\end{equation}
 We can observe that because of fourth term in the right hand side of \eqref{e51}, sparsity of the Jacobian matrix $J_1$ is lost \cite{[sk]}. In order to retain the sparsity of the Jacobian matrix, we reformulate our problem as follows \cite{[r3]}:
 Find $U^n  \in X_h$ and $d \in \mathbb{R}$ such that
 \begin{equation}\label{e52}
 l(U^n) - d \, = \, 0,
 \end{equation}
 \begin{equation}\label{e53}
 \begin{split}
 ( U^n, \, \phi_i) +& \frac{1}{d_{n,1}} \left(- d_{n,n} U^0 + \sum_{k=1}^{n-1} (d_{n,k+1} - d_{n,k}) U^{n-k}, \, \phi_i \right)\\
 &+ \frac{\Gamma {(2- \alpha)}}{d_{n,1}} \, a(d) \, (\nabla U^n, \, \nabla \phi_i) \, - \, \frac{\Gamma {(2- \alpha)}}{d_{n,1}} \, \big(f(U^n), \, \phi_i\big) \, = \,0.
 \end{split}
 \end{equation}
 For applying the Newton's method in the reformulated problem, we rewrite equations \eqref{e52} and \eqref{e53} as follows:
 \begin{equation}\label{e54}
 F_i(U^n, \, d) \, = \, 0, \quad \mbox{for} \: 1 \le i \le M+1,
 \end{equation}
 where
 \begin{equation}\label{e56}
 \begin{split}
 F_i(U^n, \, d) \, &= \, ( U^n, \, \phi_i) + \frac{1}{d_{n,1}} \left(- d_{n,n} U^0 + \sum_{k=1}^{n-1} (d_{n,k+1} - d_{n,k}) U^{n-k}, \, \phi_i \right)\\
 &+ \frac{\Gamma {(2- \alpha)}}{d_{n,1}} \, a(d) \, (\nabla U^n, \, \nabla \phi_i) \, - \, \frac{\Gamma {(2- \alpha)}}{d_{n,1}} \, \big(f(U^n), \, \phi_i\big), \, \, \, \mbox{for} \: 1 \le i \le M,
 \end{split}
 \end{equation}
 \begin{equation}\label{e57}
 F_{M+1}(U^n, \, d) \, = \, l(U^n) - d. \qquad \qquad \qquad \qquad \qquad \qquad \qquad \qquad \qquad \qquad
 \end{equation}
 Now, using Newton's method for the system of equations \eqref{e54} and \eqref{e55}, we get the following matrix equation:
 \begin{gather}\label{e58}
 J
 \begin{bmatrix}
 {\alpha^n} \\ \beta
 \end{bmatrix}
 =
 \begin{bmatrix}
 A & b\\ c & \gamma
 \end{bmatrix}
 \begin{bmatrix}
 {\alpha}^n \\ \beta
 \end{bmatrix}
 =
 \begin{bmatrix}
 \bar{F} \\ F_{M+1}
 \end{bmatrix},
 \end{gather}
  where $J$ denotes the Jacobian matrix, ${\alpha^n} \, = \, [ \alpha^n_1, \alpha^n_2,...,\, \alpha^n_M]^T,$ $\bar{F} \, = \, [F_1, F_2,..., F_M ]^T$ and entries $A \, = \, A_{M \times M},$ $b \, = b_{M \times 1}$ and $c \, = \, c_{1 \times M}$ are given below.
 \begin{equation}\label{e59}
 \begin{split}
 A_{ij} \, =& \, ( \phi_j, \, \phi_i) + \frac{\Gamma {(2- \alpha)}}{d_{n,1}} \, a(d) \, (\nabla \phi_j, \, \nabla \phi_i) - \frac{\Gamma {(2- \alpha)}}{d_{n,1}} \, \big(f'(U^n) \phi_j, \, \phi_i\big), \; \: 1 \le i, \, j \le M,\\
 b_{i1} \, =& \, \frac{\Gamma {(2- \alpha)}}{d_{n,1}} \, a'(d) \, (\nabla U^n, \, \nabla \phi_i), \quad 1 \le i \le M,\\
 c_{1j} \, =& \, \int_{\Omega} \phi_j \, dx, \quad 1 \le j \le M,\\
 \gamma \, =& \, -1.
 \end{split}
 \end{equation}
 Note that the sparsity of matrix $A$ is same as the sparsity of the Galerkin matrix corresponding to following semi-linear equation.
 \begin{equation}\label{e65}
 (w,v) + (\nabla w, \nabla v) = (f(w), v).
 \end{equation}
 So, $A$ is a sparse matrix. Hence, $J$ is also a sparse matrix \cite{[sk]}. Also, \eqref{e58} admits a unique solution \cite{[r3]}.\\

 In the following theorem we show that solution of \eqref{e49} is equivalent to the solution of \eqref{e52}-\eqref{e53}.
%\begin{figure}[h!]\label{fig1}
% 	\begin{center}
% 		%\begin{array}{cc}
% 		\includegraphics[width = 15 cm]{density.jpg}
% 		%\end{array}$
% 	\end{center}
% 	%\label{pic}
% 	\caption{\emph{ \ Density and sparsity of Jacobian matrix}}
% \end{figure}
\begin{theorem}\label{th2}
If $(U^n, \, d)$ is a solution of \eqref{e52}-\eqref{e53}, then $U^n$ is a solution of \eqref{e49}. Conversely, if $U^n$ is a solution of \eqref{e49}, then $(U^n, \, d)$ is a solution of \eqref{e52}-\eqref{e53}.
\end{theorem}	
{Proof.}
Suppose $(U^n, \, d)$ is a solution of \eqref{e52}-\eqref{e53}, then $d = l(U^n)$ and putting this in \eqref{e53}, we get
\[ \begin{split}
( U^n, \, \phi_i) + \frac{1}{d_{n,1}} &\left(- d_{n,n} U^0 + \sum_{k=1}^{n-1} (d_{n,k+1} - d_{n,k}) \, U^{n-k}, \, \phi_i \right)\\
&+ \, \frac{\Gamma {(2- \alpha)}}{d_{n,1}} \, a\big(l(U^n)\big) \, (\nabla U^n, \, \nabla \phi_i)- \, \frac{\Gamma {(2- \alpha)}}{d_{n,1}} \, \big(f(U^n), \, \phi_i\big) = 0.
\end{split} \]
Hence, $U^n$ is the solution of \eqref{e49}.
The converse is obvious. \hfill  $\square$

%%%%%%%%%%%%%%%%%%%%%%%%%%%%%%%%%%

\subsection{A priori bound}
In this section we provide \emph{a priori} bound for the fully discrete solution $U^n.$ First, we write the following coercivity property of L1 scheme.
%With these notations, equation \eqref{e13} can be written as: For each $1 \le n \le N$, find $U^n \in X_h$ such that
%\begin{equation}\label{e107}
%\begin{split}
%\left( {D}^{\alpha}_{N}U^n, v_h \right) \, - \, a\big(l(U^n)\big) \, (\Delta _hU^n, \, v_h) \, =& \, \big(P_h f(U^n), \, v_h\big), \quad  \forall v_h \in X_h, \\
%(U^0, v_h) \, =& \, (u_h^0, v_h), \quad  \forall v_h \in X_h.
%\end{split}
%\end{equation}
%Since ${D}^{\alpha}_{N}U^n$, $\Delta _hU^n$ and $P_hf^n$ lie in $X_h$, equation \eqref{e107} takes the form: For each $1 \le n \le N$, find $U^n \in X_h$ such that
%\begin{equation}\label{e108}
%\begin{split}
%{D}^{\alpha}_{N}U^n \, - \, a\big(l(U^n)\big) \, \Delta _hU^n \, =& \, P_h f(U^n), \quad  \forall v_h \in X_h, \\
%U^0 \, =& \, u_h^0.
%\end{split}
%\end{equation}
\begin{lemma}\cite{[ach1]} \label{l1}
	Let the functions $v^n = v( \cdot , \, t_n)$ be in $L^2(\Omega)$, for $n = 0, 1, . . . , N$. Then, the discrete L1 scheme satisfies
	\begin{equation}\label{ee12}
	  \left( {D}^{\alpha}_{N} v^n, \, v^n \right) \, \ge \, \frac{1}{2} \, {D}^{\alpha}_{N} \|v^n\|^2, \quad \mbox{for} \; \; n=1, 2, ..., N.
	\end{equation}
\end{lemma}
\vspace{5 pt}
For deriving \emph{a priori} estimates and \emph{a priori} error estimates for fully-discrete solution $U^n$, we need following discrete fractional Gr$\ddot{{o}}$nwall inequality.
%\begin{lemma}\label{l4}
%	\cite{[r5]} Let $(g^n)^N_{n=1}$ and $(\lambda _l)^{N-1}_{l=1}$ be given non-negative sequences. Assume that there exists a constant $\Lambda$ (independent of the step sizes) such that $ \Lambda \ge \sum_{l=0}^{N-1} \lambda_l$, and that the maximum step size satisfies
%	\begin{equation}\label{ee8}
%	\max_{1 \le n \le N } \tau_n \le \frac{1}{\sqrt[\alpha]{2 \Gamma (2- \alpha) \Lambda}}.
%	\end{equation}
%	Then, for any non-negative sequence $(v^k)^N_{k=0}$ such that
%	\begin{equation}\label{ee9}
%	D_N^{\alpha} (v^n)^2 \, \le \sum_{k=1}^{n} \lambda_{n-k} \, (v^k)^2 \, + \, v^n g^n, \quad \mbox{for} \; \, 1 \le n \le N,
%	\end{equation}
%	it holds that
%	\begin{equation}\label{ee10}
%	\begin{split}
%	v^n \le 2 E_{\alpha}(2 \,\Lambda \, t^{\alpha}_n) \, \Big( v^0  + \Gamma (1- \alpha) &\max_{1 \le j \le n} (t^{\alpha}_j g^j )\Big), \quad \mbox{for} \; \, 1 \le n \le N.
%	\end{split}
%	\end{equation}
%\end{lemma}
%\begin{remark}\label{rmk1}
%	\cite{[r5]} Above theorem remains valid if the condition \eqref{ee9} is replaced by
%	\begin{equation}\label{ee11}
%	D_N^{\alpha} (v^n) \, \le \sum_{k=1}^{n} \lambda_{n-k} \, v^k + g^n, \quad \mbox{for} \; \, 1 \le n \le N.
%	\end{equation}
%\end{remark}
\begin{lemma}\label{l5}
	\cite{[ach1], [r5a]} Let $(\xi ^n)^N_{n=1}$, $(\eta ^n)^N_{n=1}$ and $(\lambda _l)^{N-1}_{l=1}$ be given non-negative sequences. Assume that there exists a constant $\Lambda$ (independent of the step sizes) such that $ \Lambda \ge \sum_{l=0}^{N-1} \lambda_l$, and the maximum step size satisfies
	\begin{equation}\label{eee8}
	\max_{1 \le n \le N } \tau_n \le \frac{1}{\sqrt[\alpha]{2 \Gamma (2- \alpha) \Lambda}}.
	\end{equation}
    Then, for any non-negative sequence $(v^k)^N_{k=0}$ such that
	\begin{equation}\label{e105}
	\begin{split}
	  D_N^{\alpha} (v^n)^2 \, \le \, \sum_{s=1}^{n} \, \lambda_{n-s} \, (v^s)^2 + \, \xi^n \, v^n \, + \,  (\eta ^n)^2, \quad \mbox{for} \; \, 1 \le n \le N,
	\end{split}
	\end{equation}	
	it holds that $\mbox{for} \; \, 1 \le n \le N,$
	\begin{equation}\label{e109}
	\begin{split}
	  v^n \, \le \, 2 E_{\alpha} (2 \, \Lambda \, t_n^{\alpha}) \, \left[ v^0 + \max_{1 \le j \le n} \, \sum_{s=1}^{j} \, p^{(j)}_{j-s} \, \xi ^s + \sqrt{\Gamma {(1- \alpha)}} \, \max_{1 \le s \le n} \, \left\lbrace t_s^{\frac{\alpha}{2}} \, \eta ^s \right\rbrace \right].
	\end{split}
	\end{equation}	
\end{lemma}

\begin{theorem}\label{th4}
	Let $U^n$ be the solution of \eqref{e13}, then $U^n$ satisfies following estimates.
	\begin{equation}
	   \| {U^n} \| \, \le \,  C \, \big(\| {U^0} \| + 1\big),
	\end{equation}
	\begin{equation}
       \| {\nabla U^n} \| \, \le \, C \, \big( \| {\nabla U^0} \| + 1\big),
	\end{equation}
	where $n=1, 2, ..., N.$
\end{theorem}
{Proof.} Putting $v_h = U^n$ in \eqref{e13} to get
\begin{equation}\label{een7}
\begin{split}
\left( {D}^{\alpha}_{N}U^n, U^n \right) \, + \, a\big(l(U^n)\big) \, (\nabla U^n, \nabla U^n) \, =& \, \big(f(U^n), U^n\big).\\
\end{split}
\end{equation}
Using the bound of $a$ and Cauchy-Schwarz inequality in equation \eqref{een7}, we have
\begin{equation}\label{ee7}
\begin{split}
\left( {D}^{\alpha}_{N}U^n, U^n \right) +m \, ||\nabla U^n||^2 \, \le& \, \| {f(U^n)} \| \,  \| {U^n} \|.
\end{split}
\end{equation}
Lipschitz continuity of $f$ yields
%%%%%  using f(0)
%\begin{equation*}
%\begin{split}
%&\Big|\| {f(U^n)} \| - \| {f(0)} \| \Big| \, \le \| {f(U^n)-f(0)} \| \,  \le \, K \, \|{U^n} \|.
%\end{split}
%\end{equation*}
%Therefore,
%\begin{equation}\label{e17}
%\begin{split}
% \| {f(U^n)} \| \, \le \, \| {f(0)} \| + K \, \|{U^n} \| \, \le \,  C_1 \, \big(1+\|{U^n} \|\big),\\
% \end{split}
% \end{equation}
%where $C_1 = max \left\lbrace K, \|{f(0)}\| \right\rbrace.$  \\
%%%%%  using f(U^0)
\begin{equation*}
\begin{split}
&\Big|\| {f(U^n)} \| - \| {f(U^0)} \| \Big| \, \le \| {f(U^n)-f(U^0)} \| \,  \le \, K \, \big(\|{U^n} \| + \|{U^0}\|\big).
\end{split}
\end{equation*}
Therefore,
\begin{equation}\label{e17}
\begin{split}
 \| {f(U^n)} \| \, \le \, \| {f(U^0)} \| + K \,\big(\|{U^n} \| + \|{U^0}\|\big) \, \le \,  C_1 \, \big(1+\|{U^n} \|\big),\\
 \end{split}
 \end{equation}
where
\begin{equation}\label{e131}
C_1 = max \, \big\{ K, \, \|{f(U^0)}\| + K \, \|U^0\| \big\}.
\end{equation}
From \eqref{ee7} and \eqref{e17}, we get
\begin{equation}\label{ee13}
\left( {D}^{\alpha}_{N}U^n, U^n \right) \, \le \, C_1 \, \big(1+\|{U^n} \|\big) \,  \| {U^n} \|.	
\end{equation}
Also from Lemma \ref{l1}, we know that
\begin{equation*}
\left( {D}^{\alpha}_{N}U^n, U^n \right) \, \ge \, \frac{1}{2} \, D_N^{\alpha} {\| U^n\|}^2.
\end{equation*}
Thus, equation \eqref{ee13} can be written as
\begin{equation}\label{e14}
D_N^{\alpha} {\| U^n\|}^2 \, \le \, 2 C_1 \, \|{U^n} \|^2 \, + \,  2C_1 \, \|{U^n} \|.
\end{equation}
Using Lemma \ref{l5} $\big(with \ $ $v^n = \|U^n\|$, $\lambda_0 = 2C_1$, $\lambda_i=0 \, \mbox{for} \, i=1,2,...,(n-1)$, \ $\xi^n = 2C_1$ and $\eta^n = 0$$\big)$ in \eqref{e14}, we obtain
\begin{equation}\label{e90n}
\begin{split}
\|{U^n} \| \, &\le \, 2 \, E_{\alpha}(4 \, C_1 \, t^{\alpha}_n ) \: \left( \|U^0\| +  \max_{1 \le j \le n} \, \sum_{s=1}^{j} \, p^{(j)}_{j-s} \, \big(2 \, C_1\big) \right)\\
&\le \, 2 \, E_{\alpha}(4 \, C_1 \, t^{\alpha}_n ) \: \left( \|U^0\| + \Gamma {(1-\alpha)} \, \max_{1 \le j \le n} \, \Big\{ \big(2 \, C_1 \, t^{\alpha}_j \big ) \, \sum_{s=1}^{j} \, p^{(j)}_{j-s} \, \frac{t^{- \alpha}_s}{\Gamma {(1-\alpha)}} \Big\} \right).
%&\le \, 2 \, E_{\alpha}(4 \, C_1 \, t^{\alpha}_n ) \: \left( \|U^0\| +  \Gamma(1-\alpha) \, \max_{1 \le j \le n} \left\lbrace t^{\alpha}_j \, 2\, C_1 \right\rbrace \right)\\
%&\le \, 2 \, E_{\alpha}(4 \, C_1 \, t^{\alpha}_n ) \: \big( \|U^0\| + 2 \, C_1 \, \Gamma(1-\alpha) \,  T^{\alpha} \big).\\
\end{split}
\end{equation}
Using Lemma \ref{l7} in equation \eqref{e90n}, we obtain
\begin{equation}\label{e90}
\begin{split}
\|{U^n} \| \,
&\le \, 2 \, E_{\alpha}(4 \, C_1 \, t^{\alpha}_n ) \: \left( \|U^0\| +  \Gamma(1-\alpha) \, \max_{1 \le j \le n} \left\lbrace t^{\alpha}_j \, 2\, C_1 \right\rbrace \right)\\
&\le \, 2 \, E_{\alpha}(4 \, C_1 \, t^{\alpha}_n ) \: \big( \|U^0\| + 2 \, C_1 \, \Gamma(1-\alpha) \,  T^{\alpha} \big).\\
\end{split}
\end{equation}
Hence,
\begin{equation}\label{e15}
\|{U^n} \| \, \le \, C \, \big(\| {U^0} \| + 1\big),
\end{equation}
where $C = 2 \, E_{\alpha}(4 \, C_1 \, t^{\alpha}_n ) \;  max \left\lbrace 1, \, 2 \, \Gamma(1-\alpha) \, C_1 \, T^{\alpha} \right\rbrace.$\\
\\
Next, we take $v_h=\, {D}^{\alpha}_{N}U^n$ in $\eqref{e13}$ to get
\begin{equation}\label{e78}
\begin{split}
\left( {D}^{\alpha}_{N}U^n, \, {D}^{\alpha}_{N}U^n \right) \, + \, a\big(l(U^n)\big) \, (\nabla U^n, \, \nabla ( {D}^{\alpha}_{N}U^n)) \, =& \, \big(f(U^n), \, {D}^{\alpha}_{N}U^n\big).\\
\end{split}
\end{equation}
Dividing both the sides of \eqref{e78} by $a(l(U^n))$, we get
\begin{equation}\label{e78n}
\begin{split}
\frac{1}{a\big(l(U^n)\big)} \, ( {D}^{\alpha}_{N}U^n, \, {D}^{\alpha}_{N}U^n ) \, + \, \big(\nabla U^n, \, {D}^{\alpha}_{N}(\nabla U^n)\big) \, =& \,  \frac{1}{a\big(l(U^n)\big)} \, \big(f(U^n), \, {D}^{\alpha}_{N}U^n\big).\\
\end{split}
\end{equation}
 Now, in equation \eqref{e78n} we use bound of $a$ and Cauchy-Schwarz inequality to obtain
   \begin{equation}\label{e78nn}
   \begin{split}
      \frac{1}{m_2} \, \|{ {D}^{\alpha}_{N}U^n} \|^2 \, + \, \big(\nabla U^n, \, {D}^{\alpha}_{N}(\nabla U^n)\big) \, \le \, \frac{1}{m_1} \,  \|{f(U^n)} \| \, \| { {D}^{\alpha}_{N}U^n} \|. \\
   \end{split}
   \end{equation}
For $a,b>0$, using the inequality $ab \le \frac{\epsilon}{2} a^2 + \frac{1}{2 \epsilon} b^2$ $\big($with $\epsilon = m_2$$\big)$ in \eqref{e78nn}, we have
   \begin{equation*}
   \begin{split}
     \frac{1}{m_2} \, \|{ {D}^{\alpha}_{N}U^n} \|^2 \, +\, \big(\nabla U^n, \, {D}^{\alpha}_{N}(\nabla U^n)\big) \, \le& \,  \frac{m_2}{2 \, m_1^2} \, \|{f(U^n)}\|^2 + \,  \frac{1}{2 \, m_2} \, \|{ {D}^{\alpha}_{N}U^n}\|^2.\\
   \end{split}
   \end{equation*}
   This gives us
   \begin{equation}\label{e16}
   \begin{split}
       \big(\nabla U^n, \, {D}^{\alpha}_{N}(\nabla U^n)\big) \, \le& \,  \frac{m_2}{2 \, m_1^2} \, \|{f(U^n)}\|^2.\\
   \end{split}
   \end{equation}
   From equation \eqref{e17} and Poincar$\acute{e}$ inequality, we can get
   \begin{equation}
   \|{f(U^n)}\| \, \le \, C_1 \, \big(1+\|{U^n}\|\big) \, \le \, C_1 \, \big(1+ C_2\|{\nabla U^n}\|\big) \, \le \, C_3 \, \big(1+ \|{\nabla U^n}\|\big),
   \end{equation}
   where $C_2$ is a constant which appears in Poincar$\acute{e}$ inequality and $C_3= C_1 \, max \left\lbrace 1, C_2 \right\rbrace $.\\
   Also, Lemma \ref{l1} gives
   \begin{equation}\label{e16n}
      \big( \nabla U^n, {D}^{\alpha}_{N} (\nabla U^n) \big) \, \ge \, \frac{1}{2}  \, D_N^{\alpha} {\|\nabla U^n\|}^2.
   \end{equation}
From \eqref{e16} and \eqref{e16n}, we obtain
   \begin{equation}\label{e18}
   \begin{split}
   D_N^{\alpha} {\|\nabla U^n\|}^2 \, \le \, \frac{C^2_3 m_2}{m_1^2} \, \big(1+ \|{ \nabla U^n}\|\big)^2 \, \le \, C_4 \, \big(1+ {\|{ \nabla U^n}\|}^2\big),\\
   \end{split}
   \end{equation}
   where constant $C_4$ is depending on $C_3$, $m_1$ and $m_2.$\\
Using Lemma \ref{l5} (with $v^n = \| \nabla U^n\|^2$, $\lambda_0 = C_4$, $\lambda_i=0, \, \forall i=1,2,...,(n-1)$, \ $\xi^n = 0$ and $\eta^n = C_4$) in \eqref{e18} to get
  \begin{equation}\label{e91}
 \begin{split}
 \|{ \nabla U^n}\|^2 \, &\le \, 2 \, E_{\alpha}(2 \, C_4 \, t^{\alpha}_n ) \: \left( \|U^0\|^2 + \sqrt{\Gamma(1-\alpha)} \, \max_{1 \le s \le n} \left\lbrace t^{\frac{\alpha}{2}}_s \, C_4 \right\rbrace \right)\\
 &\le \, 2 \, E_{\alpha}(2 \, C_4 \, t^{\alpha}_n ) \: \Big(\|U^0\|^2 +  C_4 \, \sqrt{\Gamma(1-\alpha)} \,T^{\frac{\alpha}{2}} \Big).\\
 \end{split}
 \end{equation}
Therefore,
   \begin{equation}\label{e92}
   \begin{split}
     \|{ \nabla U^n}\|^2 \, \le  \, C \, \big( \|{ \nabla U^0}\|^2 + 1\big),
   \end{split}
   \end{equation}
   where $C = 2 \, E_{\alpha}(2 \, C_4 \, t^{\alpha}_n ) \;  max \left\lbrace 1, \, C_4 \, \sqrt{\Gamma(1-\alpha)} \,T^{\frac{\alpha}{2}} \right\rbrace.$\\
From $a^2+b^2 \, \le \, (a+b)^2$,
  \begin{equation}\label{e93}
  \begin{split}
  \|{ \nabla U^n}\|^2 \, \le  \, C \, \big( \|{ \nabla U^0}\| +1\big)^2.
  \end{split}
  \end{equation}
Hence,
   \begin{equation}\label{e19}
   \begin{split}
     \|{ \nabla U^n}\| \, \le  \, C \, \big( \|{ \nabla U^0}\| +1\big).\\
   \end{split}
   \end{equation}
This completes the proof. \hfill $\square$
\section{Existence-uniqueness of fully-discrete solution}
In this section we prove the existence and uniqueness of fully-discrete solution of the problem \eqref{e13}. For this, we use following proposition which is a consequence of Brouwer fixed point theorem \cite{[vth]}.
\begin{prop}\label{p1}
	Let $H$ be a finite dimensional Hilbert space with scalar product $(\cdot,\cdot)$ and norm $||\cdot||_H.$ Let $S:H\rightarrow H $ be a continuous map with with following properties: there exist $\rho >0$ such that
	$$(S(v),v)>0\qquad \forall v \in H \qquad \mbox{with}\qquad ||v||_H=\rho$$
	Then, there exists an element $w\in H$ such that $$S(w)=0\qquad ||w||_H\leq \rho.$$
\end{prop}

We also define
\begin{equation}\label{e130}
 \tau = \max \Bigg\{\frac{1}{\sqrt[\alpha]{C_1 \Gamma (2- \alpha)}}, \: \sqrt[\alpha]{\frac{4 \, m_1}{(K\,R_1 + L)^2 \, \Gamma (2- \alpha)}} \, \Bigg\},
\end{equation}
where $C_1$ is given in \eqref{e131}.\\
In the following  we discuss the existence and uniqueness of the fully-discrete solution.
\begin{theorem}\label{th6}
   Let $U^0, U^1,..., U^{n-1}$ are given and  $\max_{1 \le n \le N } \, \tau_n \le \tau$, then for all $1 \le n \le N$, there exists a unique solution $U^n \in X_h$ of \eqref{e13}.
\end{theorem}
{Proof.} Rewriting \eqref{e13} as follows
\begin{equation}\label{e66}
\begin{split}
( U^n, \, v_h) \, +& \, \frac{1}{d_{n,1}} \left(- d_{n,n} U^0 + \sum_{k=1}^{n-1} (d_{n,k+1} - d_{n,k}) U^{n-k}, \, v_h \right)\\
 &+ \, \frac{\Gamma {(2- \alpha )}}{d_{n,1}} \, a\big(l(U^n)\big) \, (\nabla U^n, \, \nabla v_h) - \frac{\Gamma {(2- \alpha )}}{d_{n,1}} \big(f(U^n), \, v_h\big) = \, 0.\\
\end{split}
\end{equation}
Now, we define a map $G : X_h \longrightarrow X_h$ such that
\begin{equation}\label{e67}
  \begin{split}
	\big( G(X^n), v_h \big) \, = \, &( X^n, \, v_h) + \frac{1}{d_{n,1}} \left(- d_{n,n} U^0 + \sum_{k=1}^{n-1} (d_{n,k+1} - d_{n,k}) U^{n-k}, \, v_h \right) \\
	&+ \frac{\Gamma {(2- \alpha )}}{d_{n,1}} \, a\big(l(X^n)\big) \, (\nabla X^n, \, \nabla v_h)- \frac{\Gamma {(2- \alpha )}}{d_{n,1}} \big(f(X^n), \, v_h\big).\\
  \end{split}
\end{equation}
Then the map $G$ is continuous. By choosing $v_h = X^n$ in \eqref{e67}, we get
\begin{equation}\label{e68}
\begin{split}
\big( G(X^n), X^n \big) \, = \, &( X^n, \, X^n) + \frac{1}{d_{n,1}} \left(- d_{n,n} U^0 + \sum_{k=1}^{n-1} (d_{n,k+1} - d_{n,k}) U^{n-k}, \, X^n \right)   \\
&- \frac{\Gamma {(2- \alpha )}}{d_{n,1}} \big(f(X^n), \, X^n\big) + \frac{\Gamma {(2- \alpha )}}{d_{n,1}} \, a\big(l(X^n)\big) \, (\nabla X^n, \, \nabla X^n).\\
\end{split}
\end{equation}
Applying the bound of $a$ and Cauchy-Schwarz inequality in \eqref{e68} we can arrive at
\begin{equation}\label{e69}
\begin{split}
\big( G(X^n), X^n \big) \, \ge \, \| X^n \|^2 &+ \frac{m_1 \, \Gamma {(2- \alpha )}}{d_{n,1}} \, \|\nabla X^n\|^2 - \frac{\Gamma {(2- \alpha )}}{d_{n,1}} \|f(X^n)\| \, \|X^n\| \\
&- \frac{d_{n,n}}{d_{n,1}} \, \| U^0 \| \, \|X^n\|- \frac{1}{d_{n,1}} \, \sum_{k=1}^{n-1} (d_{n,k} - d_{n,k+1})  \| U^{n-k} \| \, \|X^n \|.\\
\end{split}
\end{equation}
 Since $\frac{m_1 \, \Gamma {(2- \alpha )}}{d_{n,1}} \| \nabla X^n \|^2 \ge 0$, it follows that
\begin{equation}\label{e70}
\begin{split}
\big( G(X^n), X^n \big) \, \ge \, \Bigg(   \| X^n \| - &\frac{\Gamma {(2- \alpha )}}{d_{n,1}} \|f(X^n)\| - \frac{d_{n,n}}{d_{n,1}} \, \| U^0 \|\\
 &- \frac{1}{d_{n,1}} \, \sum_{k=1}^{n-1} (d_{n,k} - d_{n,k+1})  \| U^{n-k} \| \Bigg) \|X^n\|.\\
\end{split}
\end{equation}
From \eqref{e17} and \eqref{e70}, we get
\begin{equation}\label{e71}
  \begin{split}
    \big( G(X^n),X^n \big) \ge \Bigg\{ \left( 1- \frac{C_1 \, \Gamma {(2- \alpha )}}{d_{n,1}} \right) &\| X^n \| - \frac{C_1 \, \Gamma {(2- \alpha )}}{d_{n,1}} - \frac{d_{n,n}}{d_{n,1}} \| U^0 \|\\
     &-  \frac{1}{d_{n,1}} \sum_{k=1}^{n-1} (d_{n,k} - d_{n,k+1}) \| U^{n-k} \| \Bigg\} \|X^n\|.\\
  \end{split}
\end{equation}
Thus, $\big( G(X^n),X^n \big) > 0$ if
\begin{equation}\label{e72}
  \begin{split}
    \left( 1- \frac{C_1 \, \Gamma {(2- \alpha )}}{d_{n,1}} \right) \| X^n \| \: - \: &\frac{C_1 \, \Gamma {(2- \alpha )}}{d_{n,1}} - \frac{d_{n,n}}{d_{n,1}} \| U^0 \|\\
     &- \frac{1}{d_{n,1}} \sum_{k=1}^{n-1} (d_{n,k} - d_{n,k+1}) \| U^{n-k} \| > 0.
  \end{split}
\end{equation}
Since $\max_{1 \le n \le N } \, \tau_n \le \tau$,  $\left( 1- \frac{C_1 \, \Gamma {(2- \alpha )}}{d_{n,1}} \right) > 0.$\\
 Thus, $ \exists \; X^n \in X_h$ such that
\begin{equation}\label{e73}
\begin{split}
 \| X^n \| > \frac{1}{\left( 1- \frac{C_1 \, \Gamma {(2- \alpha )}}{d_{n,1}} \right)} \Bigg( \frac{C_1 \, \Gamma {(2- \alpha )}}{d_{n,1}} &+ \frac{d_{n,n}}{d_{n,1}} \| U^0 \| \\
 &+ \frac{1}{d_{n,1}} \sum_{k=1}^{n-1} (d_{n,k} - d_{n,k+1}) \| U^{n-k} \| \Bigg),
\end{split}
\end{equation}
Therefore, it is easy to see that $\big( G(X^n),X^n \big) > 0, \: \forall \, X^n \in X_h$ with $\|X^n\| = \rho,$\\
 where \[\rho = \frac{1}{\left( 1- \frac{C_1 \, \Gamma {(2- \alpha )}}{d_{n,1}} \right)} \left( \frac{C_1 \, \Gamma {(2- \alpha )}}{d_{n,1}} + \frac{d_{n,n}}{d_{n,1}} \| U^0 \| + \frac{1}{d_{n,1}} \sum_{k=1}^{n-1} (d_{n,k} - d_{n,k+1}) \| U^{n-k} \| \right).\]
Thus, by Proposition \ref{p1}, we can assure that \eqref{e13} has a solution.\\
Now, we prove the uniqueness of solution. For this we assume that for given $U^0, U^1, ...,$ $U^{n-1},$ there exist two solutions of \eqref{e13}, say $U^n_1$ and $U^n_2$ at time $t=t_n$. Throughout the proof, we denote $U^n_1$ by $U_1$ and $U^n_2$ by $U_2$ respectively. Let $U_1 - U_2 = r$.\\
From \eqref{e66} we can get
\begin{equation}\label{e74}
  \begin{split}
    ( U_1 - U_2, \, v_h)  + \frac{\Gamma {(2- \alpha )}}{d_{n,1}} \, a\big(l(U_1)\big) \, (\nabla U_1, \, \nabla v_h)& - \frac{\Gamma {(2- \alpha )}}{d_{n,1}} \, a(l(U_2)) \, (\nabla U_2, \, \nabla v_h)\\
    =& \, \frac{\Gamma {(2- \alpha )}}{d_{n,1}} \, \big(f(U_1)-f(U_2), \, v_h\big).\\
  \end{split}
\end{equation}
By subtracting $\frac{\Gamma {(2- \alpha )}}{d_{n,1}} \: a\big(l(U_1)\big) \, (\nabla U_2, \, \nabla v_h)$ in both sides of \eqref{e74}, we get
\begin{equation}\label{e79}
  \begin{split}
  ( r, \, v_h)  + \frac{\Gamma {(2- \alpha )}}{d_{n,1}} \, a\big(l(U_1)\big) \, (\nabla r, \, \nabla v_h) = \frac{\Gamma {(2- \alpha )}}{d_{n,1}} \, &\Big( a(l(U_2)) - a\big(l(U_1)\big)\Big)  \, (\nabla U_2, \, \nabla v_h)\\
  &+ \, \frac{\Gamma {(2- \alpha )}}{d_{n,1}} \, \big(f(U_1)-f(U_2), \, v_h\big).\\
  \end{split}
\end{equation}
Now, we take $v_h = r$ in above equation to get
\begin{equation}\label{e80}
\begin{split}
( r, \, r)  + \frac{\Gamma {(2- \alpha )}}{d_{n,1}} \, a\big(l(U_1)\big) \, (\nabla r, \, \nabla r) = \frac{\Gamma {(2- \alpha )}}{d_{n,1}} \, &\Big( a(l(U_2)) - a\big(l(U_1)\big)\Big)  \, (\nabla U_2, \, \nabla r)\\
 &+ \, \frac{\Gamma {(2- \alpha )}}{d_{n,1}} \, \big(f(U_1)-f(U_2), \, r\big).\\
\end{split}
\end{equation}
Applying the bound of $a$ and Cauchy-Schwarz inequality in \eqref{e80}, we have
\begin{equation}\label{e81}
  \begin{split}
    \|r\|^2  + \frac{m \, \Gamma {(2- \alpha )}}{d_{n,1}} \, \|\nabla r\|^2 = \frac{\Gamma {(2- \alpha )}}{d_{n,1}} \, &\big|a(l(U_2)) - a\big(l(U_1)\big)\big|  \, \|\nabla U_2\| \, \|\nabla r\|\\
    &+ \, \frac{\Gamma {(2- \alpha )}}{d_{n,1}} \, \|f(U_1)-f(U_2)\| \, \|r\|.\\
  \end{split}
\end{equation}
 Since $a$ and $f$ are Lipschitz continuous, we can get
\begin{equation}\label{e82}
  \begin{split}
     \big|a(l(U_2)) - a\big(l(U_1)\big)\big| \, \le \, L \, \|l(U_2) - l(U_1)\| \, \le \, L \, \|U_2 - U_1\| \, = \, L \, \|r\| . \\
  \end{split}
\end{equation}
\begin{equation}\label{e83}
\begin{split}
\big|f(U_1) - f(U_2)\big| \, \le \, K \, \|U_1 - U_2\| \, = \, K \, \|r\|.\\
\end{split}
\end{equation}
Also, from Theorem \ref{th4} one can get
\begin{equation}\label{e84}
\begin{split}
 \|\nabla U_2\| \, = \, R_1, \quad where \; R_1 = C(1+\|U^0\|).\\
\end{split}
\end{equation}
Therefore, from \eqref{e82}, \eqref{e83}, \eqref{e84} and Poincar$\acute{e}$ inequality, we get
\begin{equation}\label{e85}
 \begin{split}
   \|r\|^2 + \frac{m_1 \, \Gamma {(2- \alpha )}}{d_{n,1}} \, \|\nabla r\|^2 \, \le \, \left\lbrace \left( \frac{KR_1 \, \Gamma {(2- \alpha )}}{d_{n,1}} + \frac{L \, \Gamma {(2- \alpha )}}{d_{n,1}} \right) \, \|r\| \right\rbrace \, \|\nabla r\|.
 \end{split}
\end{equation}
For $a,b>0$, using the inequality $ab \le \frac{\epsilon}{2} a^2 + \frac{1}{2 \epsilon} b^2$ with $ \epsilon = \frac{d_{n,1}}{2 \, m_1 \, \Gamma {(2- \alpha )}}$ in \eqref{e85}, we get
\begin{equation}\label{e86}
 \begin{split}
   \|r\|^2  \, \le \, \frac{(KR_1 + L)^2 \, \Gamma {(2- \alpha )}}{4 \, m_1 \,d_{n,1}} \, \|r\|^2.
 \end{split}
\end{equation}
This gives us,
\begin{equation}\label{e87}
 \begin{split}
   \left( 1- \frac{(KR_1 + L)^2 \, \Gamma {(2- \alpha )}}{4 \, m_1 \, d_{n,1}} \right)  \, \|r\|^2 \, \le \, 0.
  \end{split}
\end{equation}
Since $\max_{1 \le n \le N } \, \tau_n \le \tau$, $\left( 1- \frac{(KR_1 + L)^2 \, \Gamma {(2- \alpha )}}{4 \, m_1 \, d_{n,1}} \right) > 0.$ Using this condition in \eqref{e87}, we get
\begin{equation}\label{e88}
 \begin{split}
     \|r\|^2 \, \le \, 0.
 \end{split}
\end{equation}
This shows that
\begin{equation}\label{e89}
 \begin{split}
   U_1 \, = \, U_2.
 \end{split}
\end{equation}
This completes the proof.  \hfill $\square$
\section{Error estimates}
In this section, we derive \emph{a priori} error estimate for fully-discrete solution $U^n$. For this, some additional regularity on solution $u$ is required. Therefore, we assume that there exist constants $R_2, R_3, R>0$ such that
\begin{equation}\label{e60}
\|{^{c}_{0}{D}^{\alpha}_{t_n} u}\|_{L^{\infty}({H^2(\Omega)})} \, \le \, R_3 \, , \|{\Delta u}\|_{L^{\infty}({L^2(\Omega)})} \, \le \, R_2 \quad \mbox{and} \quad \|{\nabla u}\|_{L^{\infty}({L^2(\Omega)})} \, \le \, R.
\end{equation}
For our further analysis, we recall the definition of Ritz-projection and Discrete Laplacian operators.\\
%\begin{definition}\cite{[lrs]}
%	The $L^2 - projection$ is a map  $P_h : \, L^2(\Omega) \, \rightarrow \, X_h$ defined by
%	\begin{equation}\label{e103}
%	(P_hu, \, v_h) \, = \, (u, \, v_h), \quad \forall \: v_h \, \in \, X_h.
%	\end{equation}
%\end{definition}
\begin{definition}\cite{[r6]}
	The Ritz-projection is a map  $R_h : H^1_0(\Omega) \rightarrow X_h$ such that
	\begin{equation}\label{e62}
	(\nabla w, \, \nabla v) \, = \, (\nabla R_h w, \, \nabla v), \quad \forall w \in  H^1_0(\Omega) \: \mbox{and}  \; \;  \forall v \in X_h.
	\end{equation}
\end{definition}
\noindent It is easy to prove that $R_h $  satisfies
\begin{equation}\label{e63}
\|{\nabla R_h w}\| \, \le \, R,
\end{equation}
where $R$ is given in \eqref{e60}.\\
\begin{lemma}\cite{[r7]}\label{l2}
	There exists a positive constant C (independent of h) such that
	\begin{equation}\label{e61}
	\begin{split}
	\|{w-R_hw} \|_{L^2(\Omega)} \, + h \, \|{\nabla (w-R_hw)}\|_{L^2(\Omega)} \, \le& \, Ch^2 \, \|{\Delta w}\|_{L^2(\Omega)}, \quad \forall w \in H^2 \cap H^1_0.\\
	\end{split}
	\end{equation}
\end{lemma}
\begin{definition} \cite{[vth]}
	The discrete \textit{Laplacian} is a map $\Delta _h : \, X_h \, \rightarrow \, X_h $ as
	\begin{equation}\label{e106}
	(\Delta _hu, \, v) \, = \, -(\nabla u, \, \nabla v), \quad \forall \: u, \, v \, \in \, X_h.
	\end{equation}
\end{definition}
Now, with the help of the projection operator $R_h,$ we split the error in two parts and it is given below.
\begin{equation}\label{e61n}
 	\begin{split}
u^n - U^n = u^n - R_hu^n + R_hu^n - U^n = \rho ^n + \theta ^n,
\end{split}
 	\end{equation}
 where \  $ u^n := u(t_n), \: \,  \rho ^n := u^n - R_hu^n$ \ and \ $\theta ^n := R_hu^n - U^n.$\\

	Next, in the following theorem, we provide the convergence estimate for the fully-discrete solution.

\begin{theorem}\label{th5}
	Let $u^n$ and $U^n$ be the solution of \eqref{e1} and \eqref{e13} respectively, then
	\begin{equation}\label{e20}
	\|{u^n - U^n}\| \, \le \, C \,\big(h^2 + N^{-min \left\lbrace 2- \alpha, \, r \alpha \right\rbrace }\big),
	\end{equation}
	\begin{equation}\label{e21}
	\|{ \nabla (u^n - U^n)}\| \, \le \, C \, \big(h + N^{-min \left\lbrace 2- \alpha, \, r \alpha \right\rbrace }\big),
	\end{equation}
	where $n=1, 2, ..., N.$
\end{theorem}
{Proof.}  For simplicity, we write $a := a(l(u^n))$ \ and \ $a_h := a(l(U^n)).$
For any $v_h \in X_h$ the estimate for $\theta ^n$ is given by
\begin{equation}\label{e22}
\begin{split}
&\left( {D}^{\alpha}_{N} \theta ^n, v_h \right) \, + \, a_h \, (\nabla \theta ^n, \nabla v_h) \\
 &= \, \left( {D}^{\alpha}_{N} R_hu^n, v_h \right) \, + \, a_h \, (\nabla R_hu^n, \nabla v_h) \, - \left( {D}^{\alpha}_{N} U^n, v_h \right) - a_h \, (\nabla U^n, \nabla v_h) \\
 &= \, \left( {D}^{\alpha}_{N} R_hu^n, v_h \right) \, + \, a_h \, (\nabla u^n, \nabla v_h) \, - a \, (\nabla u^n, \nabla v_h) + a \, (\nabla u^n, \nabla v_h) \, - \big(f(U^n), v_h\big)\\
&= \, \left( {D}^{\alpha}_{N} R_hu^n, v_h \right) \, + \, (a_h-a) \, (\nabla u^n, \nabla v_h) \, + \big(f(u^n), v_h\big) - ( ^{c}_{0}{D}^{\alpha}_{t_n}u, v_h) - \big(f(U^n), v_h\big)\\
&= \, \left( {D}^{\alpha}_{N} R_hu^n - \, ^{c}_{0}{D}^{\alpha}_{t_n} u , v_h \right) \, + (a_h-a) \, (\nabla R_hu^n, \nabla v_h) + \big(f(u^n)-f(U^n), v_h\big).
\end{split}
\end{equation}
We choose $v_h = \theta ^n$ in \eqref{e22} to get
\begin{equation}\label{e75}
\begin{split}
\left( {D}^{\alpha}_{N} \theta ^n, \theta  ^n \right) \, + \, a_h \, (\nabla \theta ^n, \nabla \theta  ^n) \, =& \, \left( {D}^{\alpha}_{N} R_hu^n - \, ^{c}_{0}{D}^{\alpha}_{t_n} u , \theta  ^n \right) \, +  (a_h-a) \, (\nabla R_hu^n, \nabla \theta  ^n)\\
&+ \big(f(u^n)-f(U^n), \theta  ^n\big).
\end{split}
\end{equation}
An application of Cauchy-Schwarz inequality in \eqref{e75} gives
\begin{equation}\label{e75n}
	\begin{split}
	     \left( {D}^{\alpha}_{N} \theta ^n, \theta  ^n \right) \, + \, a_h \, (\nabla \theta ^n, \nabla \theta  ^n) \, \le \, &\|{{D}^{\alpha}_{N} R_hu^n - \, ^{c}_{0}{D}^{\alpha}_{t_n} u}\| \|{ \theta ^n}\| + |a_h - a| \, \|{ \nabla R_hu^n}\| \|{ \nabla \theta ^n}\|\\
	    & + \|{f(u^n)-f(U^n)}\| \|{ \theta ^n}\|.\\
	 \end{split}
\end{equation}
By using the bound of $a$, triangle inequality, Poincar$\acute{e}$ inequality and \eqref{e63} in \eqref{e75n}, we have
\begin{equation}\label{e76}
\begin{split}
\left( {D}^{\alpha}_{N} \theta ^n, \theta  ^n \right) \, +& \, m_1 \, \| \nabla \theta ^n \|^2 \, \le \, \| {D}^{\alpha}_{N} R_hu^n - \ ^{c}_{0}{D}^{\alpha}_{t_n} R_h u \| \, \| \theta ^n \| + R \, |a_h - a| \, \|{ \nabla \theta ^n}\| \\
&+ C_2 \, \| ^{c}_{0}{D}^{\alpha}_{t_n} R_h u - \ ^{c}_{0}{D}^{\alpha}_{t_n} u \| \, \|{ \nabla \theta ^n}\|  + C_2 \,  \|{f(u^n)-f(U^n)}\| \, \|{ \nabla \theta ^n}\|, \\
\end{split}
\end{equation}
where $C_2$ is a constant which appears in Poincar$\acute{e}$ inequality.\\
For $a,b>0,$ using $ab \le \frac{\epsilon}{2} a^2 + \frac{1}{2 \epsilon} b^2$ (with $\epsilon = \frac{m_1}{3}$) in \eqref{e76} to get
\begin{equation}\label{e23}
\begin{split}
\left( {D}^{\alpha}_{N} \theta ^n, \theta  ^n \right) \, +& \, m_1 \, (\nabla \theta ^n, \nabla \theta  ^n) \, \le \, \| {D}^{\alpha}_{N} R_hu^n - \ ^{c}_{0}{D}^{\alpha}_{t_n} R_h u \| \, \| \theta ^n \| + \frac{m_1}{2} \, \|{ \nabla \theta ^n}\|^2\\
&+ \frac{3 \, C_5}{2 \, m_1} \, \Big( \| ^{c}_{0}{D}^{\alpha}_{t_n} R_h u - \ ^{c}_{0}{D}^{\alpha}_{t_n} u \|^2  + |a_h - a|^2 + \|{f(u^n)-f(U^n)}\|^2 \Big), \\
\end{split}
\end{equation}
where $C_5 = max\left\lbrace C_2^2, \, R^2 \right\rbrace.$\\	
From \eqref{e23}, we have
\begin{equation}\label{e24}
\begin{split}
\left( {D}^{\alpha}_{N} \theta ^n, \theta  ^n \right) \, \le \, &\| {D}^{\alpha}_{N} R_hu^n - \ ^{c}_{0}{D}^{\alpha}_{t_n} R_h u \| \, \| \theta ^n \| +  C_6 \Big( \| ^{c}_{0}{D}^{\alpha}_{t_n} R_h u - \ ^{c}_{0}{D}^{\alpha}_{t_n} u \|^2 \\
 &  \qquad \qquad \qquad \qquad \qquad \qquad  + |a_h-a|^2 + \|{f(u^n)-f(U^n)}\|^2 \Big),\\
\end{split}
\end{equation}	
where $C_6$ is a constant depending on $C_5$ and $m_1$.\\
Note that
\begin{equation}\label{e25}
\begin{split}
|a_h-a| \, =& \, |a_h(l(U^n)) - a(l(u^n))| \, \le \, L |l(U^n) - l(u^n)| \, \le \, L \, \|{U^n - u^n}\| \\
\le& \, L \, \big( \|{\rho^n}\| + \|{\theta^n}\| \big).
\end{split}
\end{equation}	
and
\begin{equation}\label{e26}
\begin{split}
\|{f(u^n) - f(U^n)}\| \, \le& \, K \, \|{ u^n - U^n }\| \, \le \, K \big( \|{\rho^n}\| + \|{\theta^n}\| \big).
\end{split}
\end{equation}
Also from \eqref{ee99} and \eqref{e61} one can get
\begin{equation}\label{e27}
\begin{split}
\| {D}^{\alpha}_{N} R_hu^n - \ ^{c}_{0}{D}^{\alpha}_{t_n} R_h u \| \, \le \, C \, n^{-min \left\lbrace 2- \alpha, \, r \alpha \right\rbrace},\\
\end{split}
\end{equation}
\begin{equation}\label{ee27}
\begin{split}
 \|{^{c}_{0}{D}^{\alpha}_{t_n} R_hu - \, ^{c}_{0}{D}^{\alpha}_{t_n} u}\| \, &\le \, C  \, h^2 \, \| \Delta \ ^{c}_{0}{D}^{\alpha}_{t_n} u \| \, \le \, C_7 \, h^2,
\end{split}
\end{equation}
where constant $C_7$ is depending on $R_3.$ Also, from Lemma \ref{l1}, we have
\begin{equation}\label{e28}
\begin{split}
\left( {D}^{\alpha}_{N} \theta ^n, \theta  ^n \right) \, \ge& \, \frac{1}{2}  \, {D}^{\alpha}_{N} \| \theta ^n \|^2.
\end{split}
\end{equation}
Using the values from \eqref{e25} - \eqref{e28} in \eqref{e24} to get
\begin{equation}\label{e110}
 \begin{split}
   {D}^{\alpha}_{N} \| \theta ^n \|^2 \, &\le \, 2 \, C \, n^{-min \left\lbrace 2- \alpha, \, r \sigma \right\rbrace} \, \|{\theta^n}\| + 2 \, C_6 \Big[ (C_7 \, h^2)^2 + (L^2 + K^2) \, \big(\|{\rho^n}\| + \|{\theta^n}\|\big)^2 \Big]\\
   &\le \, 2 \, C \, n^{-min \left\lbrace 2- \alpha, \, r \sigma \right\rbrace} \, \|{\theta^n}\| + 2 \, C_6 \, (C_7 \, h^2)^2 + 2 \, C_6 \, (L^2 + K^2) \, \|{\rho^n}\|^2 \\
   & \qquad \qquad \qquad \qquad \qquad \qquad \qquad \qquad \qquad + 2 \, C_6 \, (L^2 + K^2) \, \|{\theta^n}\|^2.  \\
 \end{split}
\end{equation}
Using \eqref{e61} in \eqref{e110} to get
\begin{equation}\label{e111}
\begin{split}
{D}^{\alpha}_{N} \| \theta ^n \|^2 \, &\le \, C_8 \, \|{\theta^n}\|^2 + 2 \, C \, n^{-min \left\lbrace 2- \alpha, \, r \alpha \right\rbrace} \, \|{\theta^n}\| + (C_9 \, h^2)^2, \\
\end{split}
\end{equation}
where constant $C_8$ is depending on $C_6$, $L$, $K$ and constant $C_9$ is depending on $C_6$, $C_7$, $L$, $K$. \\
From Lemma \ref{l5} $\big(with $ $v^n = \| \theta ^n\|$, $\lambda_0 = C_8$, $\lambda_i=0, \, \forall i=1,2,...,(n-1)$,  $\eta^n = C_9 \, h^2$, $\xi^n = 2\, C \, n^{-min \left\lbrace 2- \alpha, \, r \alpha \right\rbrace }$ $\big)$, we get
\begin{equation}\label{e94}
\begin{split}
\|{\theta ^n}\| \, &\le \, 2 \, E_{\alpha}(2 \, C_8 \, t^{\alpha}_n ) \: \bigg( \| \theta ^0\| + 2 \, C \, \max_{1 \le j \le n} \, \sum_{s=1}^{j} \, p_{j-s}^{(j)} \, s^{-min \left\lbrace 2- \alpha, \, r \alpha \right\rbrace } + \sqrt{\Gamma(1-\alpha)}\\
& \qquad \qquad \qquad \qquad \qquad \qquad \qquad \qquad \qquad \qquad \qquad \qquad \max_{1 \le s \le n} \left\lbrace t^{\frac{\alpha}{2}}_s \, C_9 \, h^2 \right\rbrace \bigg).\\
\end{split}
\end{equation}
Using Lemma \ref{l3} in \eqref{e94} to obtain
\begin{equation}\label{e95}
\begin{split}
\|{\theta ^n}\| \, &\le \, 2 \, E_{\alpha}(2 \, C_8 \, t^{\alpha}_n ) \: \Big(  \| \theta ^0\| + 2 \, C \, N^{-min \left\lbrace 2- \alpha, \, r \alpha \right\rbrace } + C_9 \, \sqrt{\Gamma(1-\alpha)} \, T^{\frac{\alpha}{2}} \, h^2 \Big). \\
\end{split}
\end{equation}
Choosing $U^0 = R_h u^0$, we get $\| \theta ^0 \| = 0$. Hence, from \eqref{e95}, we have
\begin{equation}\label{e30}
\begin{split}
 \|{\theta ^n}\| \, \le& \, C_{10} \, \big( h^2 + {N}^{-min \left\lbrace 2- \alpha, \, r \alpha \right\rbrace }\big),
\end{split}
\end{equation}
where $C_{10}= 2 \, E_{\alpha}(2 \, C_8 \, t^{\alpha}_n ) \; max \left\lbrace 2\, C, \, C_9 \, \sqrt{\Gamma(1-\alpha)} \, T^{\frac{\alpha}{2}} \right\rbrace$.\\
Thus,
\begin{equation}\label{e31}
\begin{split}
\|{u^n - U^n}\| \, \le& \, \|{\rho^n}\| + \|{\theta^n}\|\\
\le& \, C \, h^2 + C_{10} \, \big( h^2 + {N}^{-min \left\lbrace 2- \alpha, \, r \alpha \right\rbrace } \big)\\
\le& \, C \, \big( h^2 + {N}^{-min \left\lbrace 2- \alpha, \, r \alpha \right\rbrace } \big).
\end{split}
\end{equation}
Now, we will derive the error estimate in $H_0^1$-norm.
Using the definition of $\Delta_h$ and $R_h$, we can rewrite the equation \eqref{e22} as follows:
\begin{equation}\label{e77}
\begin{split}
\left( {D}^{\alpha}_{N} \theta ^n, \, v_h \right) \, - \, a_h \, (\Delta_h \theta ^n, \, v_h) \, =& \, \left( {D}^{\alpha}_{N} R_hu^n - \, ^{c}_{0}{D}^{\alpha}_{t_n} u , \, v_h \right) \, + (a_h-a) \, (\nabla u^n, \, \nabla v_h)\\
&+ \big(f(u^n)-f(U^n), \, v_h\big).
\end{split}
\end{equation}
From Green's theorem, it follows that $(\nabla u^n, \, \nabla v_h) = - (\Delta u^n, \, v_h)$. Therefore, equation \eqref{e77} takes the following form
\begin{equation}\label{sk1}
\begin{split}
\left( {D}^{\alpha}_{N} \theta ^n, \, v_h \right) \, - \, a_h \, (\Delta_h \theta ^n, \, v_h) \, =& \, \left( {D}^{\alpha}_{N} R_hu^n - \, ^{c}_{0}{D}^{\alpha}_{t_n} u , \, v_h \right) \, - (a_h-a) \, (\Delta u^n, \, v_h)\\
&+ \big(f(u^n)-f(U^n), \, v_h\big).
\end{split}
\end{equation}
Now, we take $v_h = - \Delta _h \theta ^n$ in \eqref{sk1} to get
\begin{equation}\label{e112}
\begin{split}
- \left( {D}^{\alpha}_{N} \theta ^n, \, \Delta _h \theta ^n \right) \, + \, a_h \, (\Delta_h \theta ^n, \, \Delta _h \theta ^n) \, =& \, - \left( {D}^{\alpha}_{N} R_hu^n - \, ^{c}_{0}{D}^{\alpha}_{t_n} u , \, \Delta _h \theta ^n \right) \, + (a_h-a) \\
&(\Delta u^n, \,\Delta _h \theta ^n) - \big(f(u^n)-f(U^n), \, \Delta _h \theta ^n\big).
\end{split}
\end{equation}
Using the definition of $\Delta _h$, we get
\begin{equation}\label{e113}
\begin{split}
 \left( \nabla {D}^{\alpha}_{N} \theta ^n, \, \nabla \theta ^n \right) \, + \, a_h \, \| \Delta_h \theta ^n\|^2 \, =& \, \big( \nabla ( {D}^{\alpha}_{N} R_hu^n - \, ^{c}_{0}{D}^{\alpha}_{t_n} u ), \, \nabla \theta ^n \big) \, + (a_h-a) \\
&(\Delta u^n, \,\Delta _h \theta ^n) - \big(f(u^n)-f(U^n), \, \Delta _h \theta ^n\big).
\end{split}
\end{equation}
Applying the bound of $a$, triangle inequality and Cauchy-Schwarz inequality in \eqref{e113} to get
\begin{equation}\label{sk2}
 \begin{split}
      \left( \nabla {D}^{\alpha}_{N} \theta ^n, \, \nabla \theta ^n \right) \, + \, &m_1 \, \| \Delta_h \theta ^n\|^2 \,  \le \, \| \nabla ({D}^{\alpha}_{N} R_hu^n - ^c_0 {D}^{\alpha}_{t_n} R_h u )\| \, \| \nabla \theta ^n \| \\
      &+ \| \nabla ( ^c_0 {D}^{\alpha}_{t_n} R_h u - ^{c}_{0}{D}^{\alpha}_{t_n} u)\| \, \| \nabla \theta ^n \| + |a-a_h| \, \| \Delta u^n\| \, \| \Delta_h \theta ^n \| \\
      &+ \|{f(u^n)-f(U^n)}\| \, \| \Delta_h \theta ^n\|.\\
 \end{split}
\end{equation} 	
By using \eqref{e60} and  $ab \le \frac{\epsilon}{2} a^2 + \frac{1}{2 \epsilon} b^2$ $\big($with $\epsilon = 1$ in $1^{st}$ term and $\epsilon = \frac{m_1}{2}$ in $3^{rd}$ \& $4^{th}$ term$\big)$ in \eqref{sk2}, we get
\begin{equation}\label{e114}
  \begin{split}
     &\left( \nabla {D}^{\alpha}_{N} \theta ^n, \, \nabla \theta ^n \right) \, + \, m_1 \, \| \Delta_h \theta ^n\|^2 \,  \le \, \| \nabla ({D}^{\alpha}_{N} R_hu^n - ^c_0 {D}^{\alpha}_{t_n} R_h u )\| \, \| \nabla \theta ^n \| + \frac{1}{2} \| \nabla \theta ^n \|^2 \\
     &+ \frac{1}{2} \| \nabla ( ^c_0 {D}^{\alpha}_{t_n} R_h u - ^{c}_{0}{D}^{\alpha}_{t_n} u)\|^2  +  \frac{R_2 ^2}{m_1} \, |a-a_h|^2 + \frac{1}{m_1} \, \|{f(u^n)-f(U^n)}\|^2  + \frac{m_1}{2} \, \| \Delta_h \theta ^n \|^2.\\
  \end{split}
\end{equation}
Therefore,
\begin{equation}\label{e32}
  \begin{split}
       ( {D}^{\alpha}_{N} ( \nabla \theta ^n), \nabla \theta ^n) \, \le \, &\frac{1}{2} \| \nabla \theta ^n \|^2 + \| \nabla ({D}^{\alpha}_{N} R_hu^n - ^c_0 {D}^{\alpha}_{t_n} R_h u )\| \, \| \nabla \theta ^n \| + \frac{R_2 ^2}{m_1} \, |a-a_h|^2 \\
       &+\frac{1}{2} \| \nabla ( ^c_0 {D}^{\alpha}_{t_n} R_h u - ^{c}_{0}{D}^{\alpha}_{t_n} u)\|^2  + \frac{1}{m_1} \, \|{f(u^n)-f(U^n)}\|^2.\\
  \end{split}
\end{equation}
Note that from \eqref{ee99}, we get
\begin{equation}\label{e115}
 \| \nabla ({D}^{\alpha}_{N} R_hu^n - ^c_0 {D}^{\alpha}_{t_n} R_h u )\| \, \le \, C \, n^{-min \left\lbrace 2- \alpha, \, r \alpha \right\rbrace }.
\end{equation}
Lipschitz continuity of $a$ and $f$ gives
\begin{equation}\label{e116}
\begin{split}
|a_h-a| \, =& \, |a_h(l(U^n)) - a(l(u^n))| \, \le \, L |l(U^n) - l(u^n)| \, \le \, L \, \|{U^n - u^n}\| \\
\le& \, L \, \big( \|{\rho^n}\| + \|{\theta^n}\| \big)\\
\le& \,C_2 \, L \, \big( \|{\nabla \rho^n}\| + \|{\nabla \theta^n}\| \big).
\end{split}
\end{equation}	
\begin{equation}\label{e117}
\begin{split}
\|{f(u^n) - f(U^n)}\| \, \le& \, K \, \|{ u^n - U^n }\| \, \le \, K \big( \|{\rho^n}\| + \|{\theta^n}\| \big) \, \le \, C_2 \, K \big( \|{\nabla \rho^n}\| + \|{\nabla \theta^n}\| \big).
\end{split}
\end{equation}
From \eqref{e61} one can get
\begin{equation}\label{e118}
\begin{split}
\|\nabla ({^{c}_{0}{D}^{\alpha}_{t_n} R_hu - \, ^{c}_{0}{D}^{\alpha}_{t_n} u})\| \, &\le \, C \, h \, \| \Delta \ ^{c}_{0}{D}^{\alpha}_{t_n} u \| \, \le \, C_{11} \, h,
\end{split}
\end{equation}
where constant $C_{11}$ is depending on $R_3.$
Also, from Lemma \ref{l1}, we have
\begin{equation}\label{e119}
\begin{split}
\left( {D}^{\alpha}_{N} (\nabla \theta ^n), \theta  ^n \right) \, \ge& \, \frac{1}{2}  \, {D}^{\alpha}_{N} \| \nabla \theta ^n \|^2.
\end{split}
\end{equation}
Using the values from \eqref{e115} - \eqref{e119} in \eqref{e32} to get
\begin{equation}\label{e35}
 \begin{split}
	   {D}^{\alpha}_{N} \| \nabla \theta ^n \|^2 \, \le \, &\| \nabla \theta ^n \|^2 + 2 \, C \, n^{-min \left\lbrace 2- \alpha, \, r \alpha \right\rbrace } \, \| \nabla \theta ^n \| + (C_{11} \, h)^2 + \left( \frac{C_2^2 \, L^2 R_2 ^2}{m_1} + \frac{C_2^2 \, K^2}{m_1} \right) \\
	   &  \big(\|\nabla \rho^n\| + \|\nabla \theta^n\|\big)^2.\\
  \end{split}
\end{equation}
Since $\|\nabla \rho^n\|, \, \|\nabla \theta^n\| \ge 0$, we can find a constant $C_{12}$ such that $\big(\|\nabla \rho^n\| + \|\nabla \theta^n\|\big)^2 \le C_{12} \, \big(\|\nabla \rho^n\|^2 + \|\nabla \theta^n\|^2\big)$ and from \eqref{e61}, $\|\nabla \rho^n\| \le C h$. Therefore, from \eqref{e35}, we get
\begin{equation}\label{e120}
\begin{split}
{D}^{\alpha}_{N} \| \nabla \theta ^n \|^2 \, \le \, & C_{13} \, \| \nabla \theta ^n \|^2 + 2 \, C \, n^{-min \left\lbrace 2- \alpha, \, r \alpha \right\rbrace } \, \| \nabla \theta ^n \| + (C_{14} \, h)^2.\\
\end{split}
\end{equation}
where constant $C_{13}$ is depending on $L$, $K$, $R_2$, $m_1$, $C_{2}$ and constant $C_{14}$ is depending on $L$, $K$, $R_2$, $m_1$, $C_{2}$ and $C_{11}$.\\
An application of Lemma \ref{l5} $\big($with $v^n = \| \nabla \theta ^n\|$, $\lambda_0 = C_{13}$, $\lambda_i=0, \, \forall i=1,2,...,(n-1)$,  $\xi^n = 2\, C \, n^{-min \left\lbrace 2- \alpha, \, r \alpha \right\rbrace }$, $\eta^n = C_{14} \, h$$\big)$ gives
\begin{equation}\label{e121}
\begin{split}
\|{\nabla \theta ^n}\| \, &\le \, 2 \, E_{\alpha}(2 \, C_{13} \, t^{\frac{\alpha}{2}}_n ) \: \bigg( \| \nabla \theta ^0\| + 2 \, C \, \max_{1 \le j \le n} \, \sum_{s=1}^{j} \, p_{j-s}^{(j)} \, s^{-min \left\lbrace 2- \alpha, \, r \alpha \right\rbrace } + \sqrt{\Gamma(1-\alpha)}\\
& \qquad \qquad \qquad \qquad \qquad \qquad \qquad \qquad \qquad \qquad \qquad \qquad \max_{1 \le s \le n} \left\lbrace t^{\frac{\alpha}{2}}_s \, C_{14} \, h \right\rbrace \bigg).\\
\end{split}
\end{equation}
From Lemma \ref{l3}, we get
\begin{equation}\label{e122}
\begin{split}
\|{\nabla \theta ^n}\| \, &\le \, 2 \, E_{\alpha}(2 \, C_{13} \, t^{\alpha}_n ) \: \left(  \| \nabla \theta ^0\| + 2 \, C \, N^{-min \left\lbrace 2- \alpha, \, r \alpha \right\rbrace } + C_{14} \, \sqrt{\Gamma(1-\alpha)} \, T^{\frac{\alpha}{2}} \, h \right). \\
\end{split}
\end{equation}
Choosing $U^0 = R_h u^0$, to get $\| \nabla \theta ^0 \| = 0$. Hence, from \eqref{e122}, we have
\begin{equation}\label{e123}
\begin{split}
\|{\theta ^n}\| \, \le& \, C_{15} \, \big( h + {N}^{-min \left\lbrace 2- \alpha, \, r \alpha \right\rbrace }\big),
\end{split}
\end{equation}
where $C_{15}= 2 \, E_{\alpha}(2 \, C_{13} \, t^{\alpha}_n ) \; max \left\lbrace 2\, C, \, C_{14} \, \sqrt{\Gamma(1-\alpha)} \, T^{\frac{\alpha}{2}} \right\rbrace$.\\
Hence,
\begin{equation}\label{e38}
  \begin{split}
      \|{ \nabla (u^n - U^n)}\| \, \le& \, \|{ \nabla \rho^n}\| + \|{ \nabla \theta^n}\|\\
      \le& \, C \, h + C_{15} \, \big(h + {N}^{-min \left\lbrace 2- \alpha, \, r \alpha \right\rbrace } \big)\\
      \le& \, C \, \big( h + {N}^{-min \left\lbrace 2- \alpha, \, r \alpha \right\rbrace } \big).
  \end{split}
\end{equation}
This completes the proof.  \hfill $\square$
%%%%%%%%%%%%%%%%%%%%%%%%%%%%%
\begin{corollary}
The L1-FEM solution $U^n$ satisfies:
\begin{equation}
\|{  u^n - U^n}\|_1 \leq C \, \big( h + {N}^{-min \left\lbrace 2- \alpha, \, r \alpha \right\rbrace }  \big) \ \mbox{for} \ n=0,1,..., N.
\end{equation}
\end{corollary}
{Proof.} By using Poincar$\acute{e}$ inequality, we can write
\begin{equation}
\|{  u^n - U^n}\|_1 \leq C \,\|{ \nabla (u^n - U^n)}\|.
\end{equation}
Now using the estimate \eqref{e21} from Theorem \ref{th5}, we get the desired result. \hfill $\square$

\begin{remark}
 For the time fractional partial differential equation  \eqref{cuc3:1.1}- \eqref{cuc3:1.3}, we have assumed certain regularity of solution $u$ in Lemma \ref{l6} and in equation \eqref{e60}. The derivation of assumed regularity of solution $u$ to our problem is still open.
\end{remark}

\section{Numerical experiments}
In this section, we present three different numerical examples to conform our theoretical estimates. For calculating the error, we consider problems with known exact solutions which also satisfy the assumption of Lemma \ref{l6}. In every example, we consider the time interval $[0,1]$ and tolerance $\epsilon = 10^{-12}$ for stopping the iterations in Newton's method. Let $N$ denote the number of sub-intervals in time direction and $\big(M_s+1\big)$ be the number of node points in each spatial direction. \medskip \\
{Example 1:}  Consider the time fractional non-local PDE
\begin{equation}\label{E1}
  \begin{split}
     ^{c}_{0}D^{\alpha}_{t}u(x,t) \, - \, a\big(l(u)\big) \: \Delta u(x,t) \, =& \, f(x,t,u),\;\; \mbox{in} \: \Omega \times (0,1],\\
     u(x,t)=& \, 0 \;\; \mbox{on} \; \partial  \Omega, \; \forall t \in [0,1], \\
     u(x,0)=& \, 0 \;\; \mbox{in} \; \Omega,
  \end{split}
\end{equation}
where $\Omega = (0, \pi)$, $a(x)=3+\sin x$ and we choose $f(x,t,u)$ in such a way that the exact solution of given PDE be $u(x,t)=t^3 \sin x$.\\
To obtain the order of convergence in time, we fix $M_s=1000$ and calculate the error for different values of $N.$ Table 1 shows order of convergence in the temporal direction in $L^{\infty}$ norm using uniform mesh.\\
To obtain the convergence rate in spatial direction, we fix $N=15000$ and calculate the error for different values of $M_s.$ The order of convergence in spatial direction in $L^2$ and $H^1_0$ norms are given in the Table 2 and Table 3, respectively. Through this example, we have shown that if solution $u$ does not have initial singularity then there is no issue in getting optimal order of convergence in $L^{\infty}$ norm in temporal direction.
\begin{center}
  \begin{tiny}	
	\begin{table}[h!]
		\begin{center}
		\begin{tabular}{|c|c|c|c|c|c|c|}
			\hline
	   %%%  FOR COLOR FONT
%			\multirow{2}{*}{\textcolor{red}{\large $N$}} &        \multicolumn{2}{c|}{\textcolor{red}{$\alpha = 0.4$}} & \multicolumn{2}{c|}{\textcolor{red}{$\alpha = 0.5$}} & \multicolumn{2}{c|}{\textcolor{red}{$\alpha = 0.7$}}\\
%			\cline{2-7}
%			&  \textcolor{plum}{Error} & \textcolor{plum}{OC} & \textcolor{plum}{Error} & \textcolor{plum}{OC} & \textcolor{plum}{Error} & \textcolor{plum}{OC}  \\
        %%% FOR BLACK FONT	
			\multirow{2}{*}{\large $N$} & \multicolumn{2}{c|}{$\alpha = 0.4$} & \multicolumn{2}{c|}{$\alpha = 0.5$} & \multicolumn{2}{c|}{$\alpha = 0.7$}\\
			\cline{2-7}
			& Error & OC & Error & OC & Error & OC  \\
		%%%	
			\hline
			$2^6$  & 3.40E-04 & 1.5698 & 6.39E-04 & 1.4789 & 2.08E-03 & 1.2911 \\
			\cline{1-1}\cline{2-7}
			$2^7$   & 1.15E-04 & 1.5912 & 2.29E-04 & 1.4911 & 8.49E-04 & 1.2955  \\
			\cline{1-1}\cline{2-7}
			$2^8$   & 3.80E-05 & 1.6339 & 8.15E-05 & 1.5104 & 3.46E-04 & 1.300039266 \\
			\cline{1-1}\cline{2-7}
			$2^9$ & 1.23E-05 & - & 2.86E-05 & - & 1.40E-04 & -  \\ 					
			\hline
		\end{tabular}
	    \end{center}
		\caption {\emph {Error and order of convergence in $L^{\infty}$ norm in temporal direction on uniform mesh for Example 1.}}
		\label{T1}
	\end{table}
\end{tiny}
\end{center}
%%%%%%%%%
\begin{center}
	\begin{tiny}	
		\begin{table}[h!]
		\begin{center}
			\renewcommand{\arraystretch}{1.2}
			\begin{tabular}{|c|c|c|c|c|c|c|}
				\hline
		%%%  FOR COLOR FONT
%				\multirow{2}{*}{\textcolor{red}{\large $M$}} &    \multicolumn{2}{c|}{\textcolor{red}{$\alpha = 0.4$}} & \multicolumn{2}{c|}{\textcolor{red}{$\alpha = 0.5$}} & \multicolumn{2}{c|}{\textcolor{red}{$\alpha = 0.7$}}\\
%				\cline{2-7}
%				 & \textcolor{plum}{Error} & \textcolor{plum}{OC} & \textcolor{plum}{N} & \textcolor{plum}{Error} & \textcolor{plum}{OC} & \textcolor{plum}{N} & \textcolor{plum}{Error} & \textcolor{plum}{OC}  \\
	      %%%  FOR BLACK FONT
	      			\multirow{2}{*}{\large $M_s$} &    \multicolumn{2}{c|}{$\alpha = 0.4$} & \multicolumn{2}{c|}{$\alpha = 0.5$} & \multicolumn{2}{c|}{$\alpha = 0.7$}\\
	      			\cline{2-7}
	      			 & Error & OC & Error & OC  & Error & OC  \\
	      		%%%%%%	
			    \hline
				$2^4$ & 8.76E-03 & 2.0060  & 8.77E-03 & 2.0057  & 8.60E-03 & 2.0061 \\
				\hline
				$2^5$  & 2.18E-03 & 2.0026 & 2.18E-03 & 2.0034 & 2.14E-03 & 2.0091 \\
				\hline
				$2^6$ & 5.44E-04 & 2.0053  &5.44E-04 & 2.0086 & 5.32E-04 & 2.0327 \\
				\hline
				$2^7$ & 1.36E-04 & -  & 1.35E-04 & - & 1.30E-04 & -   \\
				\hline
			\end{tabular}
			\caption {\emph {Error and order of convergence  in $L^2$ norm in space for Example 1.}}
			\label{T2}
			\medskip
		\end{center}
		\end{table}
	\end{tiny}	
\end{center}
%%%%%%%%
\begin{center}
	\begin{tiny}
		\begin{table}[h!]
			\begin{center}
			\renewcommand{\arraystretch}{1.2}
			\begin{tabular}{|c|c|c|c|c|c|c|}
				\hline
			%%%  FOR COLOR FONT
				%				\multirow{2}{*}{\textcolor{red}{\large $M$}} &    \multicolumn{2}{c|}{\textcolor{red}{$\alpha = 0.4$}} & \multicolumn{2}{c|}{\textcolor{red}{$\alpha = 0.5$}} & \multicolumn{2}{c|}{\textcolor{red}{$\alpha = 0.7$}}\\
				%				\cline{2-7}
				%				&  \textcolor{plum}{Error} & \textcolor{plum}{OC}  & \textcolor{plum}{Error} & \textcolor{plum}{OC}  & \textcolor{plum}{Error} & \textcolor{plum}{OC}  \\
			%%%  FOR BLACK FONT
				\multirow{2}{*}{\large $M_s$} &    \multicolumn{2}{c|}{$\alpha = 0.4$} & \multicolumn{2}{c|}{$\alpha = 0.5$} & \multicolumn{2}{c|}{$\alpha = 0.7$}\\
				\cline{2-7}
				 & Error & OC &  Error & OC &  Error & OC  \\
				%%%%%%	
				\hline
				$2^4$  & 1.42E-01 & 0.9992 & 1.42E-01 & 0.9992 & 1.42E-01 & 0.9992	\\
				\hline
				$2^5$  & 7.10E-02 & 0.9998 & 7.10E-02 & 0.9998 & 7.10E-02 & 0.9998	\\
				\hline
				$2^6$  & 3.55E-02 & 0.9999 & 3.55E-02 & 0.9999 & 3.55E-02 & 0.9999	\\
				\hline
				$2^7$ & 1.78E-02 & - & 1.78E-02 & - & 1.78E-02 & -	\\
				\hline	
			\end{tabular}
			\caption {\emph {Error and order of convergence  in $H^{1}_0$ norm in space for Example 1.}}
			\label{T3}
			\end{center}
		\end{table}
	\end{tiny}
\end{center}
{Example 2:}  Consider the time fractional non-local PDE
\begin{equation}\label{E2}
\begin{split}
^{c}_{0}D^{\alpha}_{t}u(x,t) \, - \, a\big(l(u)\big) \: \Delta u(x,t) \, =& \, f(x,t,u),\;\; \mbox{in} \: \Omega \times (0,1],\\
u(x,t)=& \, 0 \;\; \mbox{on} \; \partial  \Omega, \; \forall t \in [0,1], \\
u(x,0)=& \, 0 \;\; \mbox{in} \; \Omega,
\end{split}
\end{equation}
where $\Omega = (0, \pi)$, $a(x)=3+\sin x$ and we choose $f(x,t,u)$ in such a way that the exact solution of given PDE be $u(x,t)=(t^3+t^{\alpha}) \sin x$.\\
To obtain the order of convergence in temporal direction, we fix $M_s=1000$ and calculate error for different values of $N.$  Table 4 shows the order of convergence in the temporal direction in $L^{\infty}$ norm on uniform mesh. In Table 5, we provide the order of convergence in the temporal direction in $L^{\infty}$ norm on graded mesh with grading parameter  $r = \frac{2-\alpha}{\alpha}.$\\
To obtain the convergence rate in spatial direction, we fix $N=15000$ and calculate the error for different values of $M_s.$ The order of convergence in spatial direction in $L^2$ and $H^1_0$ norms are given in the Table 6 and Table 7, respectively.
\begin{center}
	\begin{tiny}
		\begin{table}[h!]
            \begin{center}
			\begin{tabular}{|c|c|c|c|c|c|c|}
				\hline
			%%%  FOR COLOR FONT
%			\multirow{2}{*}{\textcolor{red}{\large $N$}} &        \multicolumn{2}{c|}{\textcolor{red}{$\alpha = 0.4$}} & \multicolumn{2}{c|}{\textcolor{red}{$\alpha = 0.5$}} & \multicolumn{2}{c|}{\textcolor{red}{$\alpha = 0.7$}}\\
%			\cline{2-7}
%			&  & \textcolor{plum}{Error} & \textcolor{plum}{OC} & \textcolor{plum}{Error} & \textcolor{plum}{OC} & \textcolor{plum}{Error} & \textcolor{plum}{OC}  \\	
			%%%  FOR BLACK FONT
			    \multirow{2}{*}{\large $N$} & \multicolumn{2}{c|}{$\alpha = 0.4$} & \multicolumn{2}{c|}{$\alpha = 0.5$} & \multicolumn{2}{c|}{$\alpha = 0.7$}\\
			    \cline{2-7}
			     & Error & OC & Error & OC & Error & OC  \\
			    %%%	
				\hline
				$2^6$ & 3.27E-02 & 0.2463  & 2.60E-02 &  0.3830 & 1.13E-02 & 0.6384 \\
				\cline{1-1}\cline{2-7}
				$2^7$ & 2.76E-02 & 0.2860 & 1.99E-02 & 0.4207  & 7.27E-03 & 0.6638  \\
				\cline{1-1}\cline{2-7}
				$2^8$   & 2.26E-02 & 0.3163 & 1.49E-02 & 0.4461 & 4.59E-03 & 0.6784 \\
				\cline{1-1}\cline{2-7}
				$2^9$ & 1.82E-02 & - & 1.09E-02 & - & 2.87E-03 & -  \\
				\hline
			\end{tabular}
		    \end{center}
			\caption {\emph {Error and order of convergence  in $L^{\infty}$ norm in temporal direction on uniform mesh for Example 2.}}
			\label{T4}
		\end{table}
	\end{tiny}
\end{center}
%%%%%
\begin{center}
	\begin{tiny}
		\begin{table}[h!]
			\begin{center}
				\begin{tabular}{|c|c|c|c|c|c|c|}
					\hline
					%%%  FOR COLOR FONT
%								\multirow{2}{*}{\textcolor{red}{\large $N$}} &        \multicolumn{2}{c|}{\textcolor{red}{$\alpha = 0.4$}} & \multicolumn{2}{c|}{\textcolor{red}{$\alpha = 0.5$}} & \multicolumn{2}{c|}{\textcolor{red}{$\alpha = 0.7$}}\\
%								\cline{2-7}
%								&  & \textcolor{plum}{Error} & \textcolor{plum}{OC} & \textcolor{plum}{Error} & \textcolor{plum}{OC} & \textcolor{plum}{Error} & \textcolor{plum}{OC}  \\	
				%%%  FOR BLACK FONT
					\multirow{2}{*}{\large $N$} & \multicolumn{2}{c|}{$\alpha = 0.4$} & \multicolumn{2}{c|}{$\alpha = 0.5$} & \multicolumn{2}{c|}{$\alpha = 0.7$}\\
					\cline{2-7}
					 & Error & OC & Error & OC & Error & OC  \\
					%%%	
				\hline
				$2^6$ & 3.77E-03 & 1.4771 & 4.21E-03 & 1.4157  & 6.10E-03 &  1.2649 \\
				\hline
				$2^7$ & 1.36E-03 & 1.5210 & 1.58E-03 & 1.4498  & 2.54E-03 & 	1.2815  \\
				\hline
				$2^8$ & 4.72E-04 & 1.5530 & 5.78E-04 & 1.4735 & 1.04E-03 &	1.2921 \\
				\hline
				$2^9$ & 1.61E-04 & - & 2.08E-04 & - & 4.27E-04 & -  \\
				\hline
			\end{tabular}
		    \end{center}
			\caption {\emph {Error and order of convergence  in $L^{\infty}$ norm in temporal direction on graded mesh for Example 2.}}
			\label{T5}
		\end{table}
	\end{tiny}
\end{center}
%%%%%%%%
\begin{center}
	\begin{tiny}
		\begin{table}[h!]
		  \begin{center}	
			\renewcommand{\arraystretch}{1.2}
			\begin{tabular}{|c|c|c|c|c|c|c|}
				\hline
			%%%  FOR COLOR FONT
%				\multirow{2}{*}{\textcolor{red}{\large $M$}} &    \multicolumn{3}{c|}{\textcolor{red}{$\alpha = 0.4$}} & \multicolumn{3}{c|}{\textcolor{red}{$\alpha = 0.5$}} & \multicolumn{3}{c|}{\textcolor{red}{$\alpha = 0.7$}}\\
%				\cline{2-10}
%				 & \textcolor{plum}{Error} & \textcolor{plum}{OC}  & \textcolor{plum}{Error} & \textcolor{plum}{OC} & \textcolor{plum}{Error} & \textcolor{plum}{OC}  \\
			%%%  FOR BLACK FONT
				\multirow{2}{*}{\large $M_s$} &    \multicolumn{2}{c|}{$\alpha = 0.4$} & \multicolumn{2}{c|}{$\alpha = 0.5$} & \multicolumn{2}{c|}{$\alpha = 0.7$}\\
				\cline{2-7}
				 & Error & OC  & Error & OC  & Error & OC  \\
				%%%%%%	
				\hline
				$2^4$  & 4.03E-03 & 2.0014 & 3.94E-03 & 2.0018 & 3.77E-03 & 2.0041 \\
				\hline
				$2^5$ & 1.01E-03 & 2.0018 & 9.85E-04 & 2.0030 & 9.40E-04 & 2.0113 \\
				\hline
				$2^6$  & 2.51E-04 & 2.0063 & 2.46E-04 & 2.011 & 2.33E-04 & 2.0434 \\
				\hline
				$2^7$ & 6.26E-05 & - & 6.10E-05 & - & 5.65E-05 & - \\
				\hline
				\end{tabular}
			\caption {\emph{Error and order of convergence  in $L^2$ norm in space for Example 2.}}
			\label{T6}
			\end{center}
		\end{table}
	\end{tiny}
\end{center}
%%%%%%%
\begin{center}
	\begin{tiny}
		\begin{table}[h!]
			\begin{center}
			\renewcommand{\arraystretch}{1.2}
            \begin{tabular}{|c|c|c|c|c|c|c|}
	        \hline
	        %%%  FOR COLOR FONT
%				\multirow{2}{*}{\textcolor{red}{\large $M$}} &    \multicolumn{3}{c|}{\textcolor{red}{$\alpha = 0.4$}} & \multicolumn{3}{c|}{\textcolor{red}{$\alpha = 0.5$}} & \multicolumn{3}{c|}{\textcolor{red}{$\alpha = 0.7$}}\\
%				\cline{2-7}
%				 & \textcolor{plum}{Error} & \textcolor{plum}{OC} &  \textcolor{plum}{Error} & \textcolor{plum}{OC} & \textcolor{plum}{Error} & \textcolor{plum}{OC}  \\
	       %%%  FOR BLACK FONT
                \multirow{2}{*}{\large $M_s$} &    \multicolumn{2}{c|}{$\alpha = 0.4$} & \multicolumn{2}{c|}{$\alpha = 0.5$} & \multicolumn{2}{c|}{$\alpha = 0.7$}\\
	            \cline{2-7}
	            &  Error & OC &  Error & OC &  Error & OC  \\
	            %%%%%%	
				\hline
				$2^4$ & 7.10E-02 & 0.9993 & 7.10E-02 & 0.9993 &  7.10E-02 & 0.9994 \\
				\hline
				$2^5$ & 3.55E-02 & 0.9998 & 3.55E-02 & 0.9998 & 3.55E-02 & 0.9998 \\
				\hline
				$2^6$ & 1.78E-02 & 0.9999 & 1.78E-02 & 0.9999 & 1.78E-02 & 0.9999 \\
				\hline
				$2^7$ & 8.88E-03 & - & 8.88E-03 & - & 8.88E-03 & - \\
				\hline
			\end{tabular}
			\caption {\emph {Error and order of convergence  in $H^1_0$ norm in space for Example 2.}}
			\label{T7}
			\end{center}
		\end{table}
	\end{tiny}
\end{center}
{Example 3:}  Consider the time fractional non-local PDE
\begin{equation}\label{E3}
\begin{split}
^{c}_{0}D^{\alpha}_{t}u(x,t) \, - \, a\big(l(u)\big) \: \Delta u(x,t) \, =& \, f(x,t,u),\;\; \mbox{in} \: \Omega \times (0,1],\\
u(x,t)=& \, 0 \;\; \mbox{on} \; \partial \Omega, \; \forall t \in [0,1], \\
u(x,0)=& \, 0 \;\; \mbox{in} \; \Omega,
\end{split}
\end{equation}
where $\Omega = (0,1) \times (0,1)$, $a(x)=3+\sin x$ and we choose $f(x,t,u)$ in such a way that the exact solution of given PDE is $u(x,t)=(t^3+t^{\alpha}) (x-x^2)(y-y^2).$\\
To obtain the order of convergence in temporal direction, we take different values of $N$ and use the relation $M_s = \big\lfloor N^{\frac{2-\alpha}{2}} \big\rfloor$ for $\alpha = 0.5$, and ,  $M_s = 2 \, \big\lfloor N^{\frac{2-\alpha}{2}} \big\rfloor$ for $\alpha = 0.7.$ In Table 8, we have given the order of convergence in the temporal direction in $L^{\infty}$ norm on graded mesh with grading parameter  $r = \frac{2-\alpha}{\alpha}.$\\
To obtain the convergence rate in spatial direction, we take $N=\big\lfloor M_s^{\frac{2}{(2-\alpha)}} \big\rfloor$. The order of convergence in spatial direction in $L^2$ and $H^1_0$ norms are given in the Table 9 and Table 10, respectively.\\
The graph of exact and numerical solutions for $\alpha = 0.5$ are shown in Figure 1.\\
%\begin{center}
%  \begin{tiny}	
%	\begin{table}[h!]\label{T8}
%      \begin{center}
%		\begin{tabular}{|c|c|c|c|c|}
%			\hline
%		%%%  FOR COLOR FONT
%%			\multirow{2}{*}{\textcolor{red}{\large $N$}} &     \multicolumn{2}{c|}{\textcolor{red}{$\alpha = 0.5$}} & \multicolumn{2}{c|}{\textcolor{red}{$\alpha = 0.7$}}\\
%%	    	\cline{2-5}
%%			 \textcolor{plum}{Error} & \textcolor{plum}{OC}  & \textcolor{plum}{Error} & \textcolor{plum}{OC}  \\
%		 %%%  FOR BLACK FONT
%		   \multirow{2}{*}{\large $N$} &     \multicolumn{2}{c|}{$\alpha = 0.5$} & \multicolumn{2}{c|}{$\alpha = 0.7$}\\
%		   \cline{2-5}
%		    & Error & OC  & Error & OC  \\	
%		  %%%	
%			\hline	
%		    $2^4$ & 3.16E-03 & 1.4909 & 1.39E-03 & 1.2815\\
%			\hline
%			$3^4$ & 2.82E-04 & 1.4995  & 1.74E-04 &  1.3042 \\
%			\hline
%			$4^4$  & 5.01E-05 & 1.500 & 3.87E-05 & 1.3250 \\
%			\hline
%			$5^4$  & 1.31E-05 & -   & 1.19E-05 &  -  \\
%			\hline
%		\end{tabular}
%	   \end{center}
%		\caption {\emph{Error and order of convergence in $L^2$ norm in time on graded mesh for Example 3.}}
%	\end{table}
%  \end{tiny}
%\end{center}
%%%%%%%
\begin{center}
  \begin{tiny}	
	\begin{table}[h!]
	  \begin{center}
	  	\begin{tabular}{|c|c|c|c|c|}
	  		\hline
	  %%%  FOR COLOR FONT
%			\multirow{2}{*}{\textcolor{red}{\large $N$}} &     \multicolumn{2}{c|}{\textcolor{red}{$\alpha = 0.5$}} & \multicolumn{2}{c|}{\textcolor{red}{$\alpha = 0.7$}}\\
%	    	\cline{2-5}
%			 & \textcolor{plum}{Error} & \textcolor{plum}{OC} &  \textcolor{plum}{Error} & \textcolor{plum}{OC}  \\
	  %%%  FOR BLACK FONT
	  	\multirow{2}{*}{\large $N$} &     \multicolumn{2}{c|}{$\alpha = 0.5$} & \multicolumn{2}{c|}{$\alpha = 0.7$}\\
	  	\cline{2-5}
	  	& Error & OC & Error & OC  \\	
	  %%%	
	  	\hline	
	    $2^4$  & 3.16E-03 & 1.4909  & 1.39E-03 & 1.2815 \\
	   \hline
	    $3^4$  & 2.82E-04 & 1.4995 & 1.74E-04 & 1.3042  \\
	   \hline
	    $4^3$ & 5.01E-05 & 1.500 &  3.87E-05 & 1.3250 \\
	   \hline
	    $5^4$   &  1.31E-05 & -  &   1.19E-05 & -  \\
	   \hline		
   \end{tabular}
   \end{center}
    \caption {\emph{Error and order of convergence in $L^{\infty}$ norm in time on graded mesh for Example 3.}}
    \label{T9}
  \end{table}
 \end{tiny}
\end{center}
%%%%%%%%%%
\begin{center}
  \begin{tiny}	
	\begin{table}[h!]
		\renewcommand{\arraystretch}{1.2}
		\begin{center}
		\begin{tabular}{|c|c|c|c|c|}
			\hline
	 %%%  FOR COLOR FONT
	%			\multirow{2}{*}{\textcolor{red}{\large $M$}} &     \multicolumn{2}{c|}{\textcolor{red}{$\alpha = 0.5$}} & \multicolumn{2}{c|}{\textcolor{red}{$\alpha = 0.7$}}\\
	%	    	\cline{2-5}
	%			 & \textcolor{plum}{Error} & \textcolor{plum}{OC} &  \textcolor{plum}{Error} & \textcolor{plum}{OC}  \\
	       %%%  FOR BLACK FONT
	         \multirow{2}{*}{\large $M_s$} &     \multicolumn{2}{c|}{$\alpha = 0.5$} & \multicolumn{2}{c|}{$\alpha = 0.7$}\\
             \cline{2-5}
	         & Error & OC & Error & OC  \\	
	      %%%	
			\hline
			$2^2$ & 1.20E-02 & 1.9223 & 1.20E-02 & 1.9192 \\
			\hline
			$2^3$ & 3.16E-03 & 1.9808 & 3.16E-03 & 1.9797 \\
			\hline
			$2^4$  & 8.00E-04 & 1.9963 & 8.01E-04 & 1.9951 \\
			\hline
			$2^5$  & 2.00E-04 & - &  2.01E-04 & - \\
			\hline	
		\end{tabular}
	    \end{center}
		\caption {\emph {Error and order of convergence in $L^2$ norm in space for Example 3.}}
		\label{T10}
	\end{table}
  \end{tiny}
\end{center}
%%%%%%%%
\begin{center}
	\begin{tiny}	
		\begin{table}[h!]
			\renewcommand{\arraystretch}{1.2}
			\begin{center}
			\begin{tabular}{|c|c|c|c|c|}
			\hline
		 %%%  FOR COLOR FONT
		%			\multirow{2}{*}{\textcolor{red}{\large $M$}} &     \multicolumn{2}{c|}{\textcolor{red}{$\alpha = 0.5$}} & \multicolumn{2}{c|}{\textcolor{red}{$\alpha = 0.7$}}\\
		%	    	\cline{2-5}
		%			 & \textcolor{plum}{Error} & \textcolor{plum}{OC} &  \textcolor{plum}{Error} & \textcolor{plum}{OC}  \\
		%%%  FOR BLACK FONT
		       \multirow{2}{*}{\large $M_s$} &     \multicolumn{2}{c|}{$\alpha = 0.5$} & \multicolumn{2}{c|}{$\alpha = 0.7$}\\
		       \cline{2-5}
		       & Error & OC & Error & OC  \\	
		     %%%	
				\hline
				$2^2$  & 8.80E-02 & 0.9802 & 8.79E-02 & 0.9801 \\
				\hline
				$2^3$ & 4.46E-02 & 0.9950 & 4.46E-02 & 0.9950 \\
				\hline
				$2^4$  & 2.24E-02 & 0.9987 & 2.24E-02 & 0.9987 \\
				\hline
				$2^5$  & 1.12E-02 & - & 1.12E-02 & - \\
				\hline
			\end{tabular}
		    \end{center}
			\caption {\emph {Error and order of convergence in $H^1_0$ norm in space for Example 3.}}
			\label{T11}
		\end{table}
	\end{tiny}
\medskip
\end{center}

\begin{figure}[h!]
	\begin{center}
		\epsfig{file=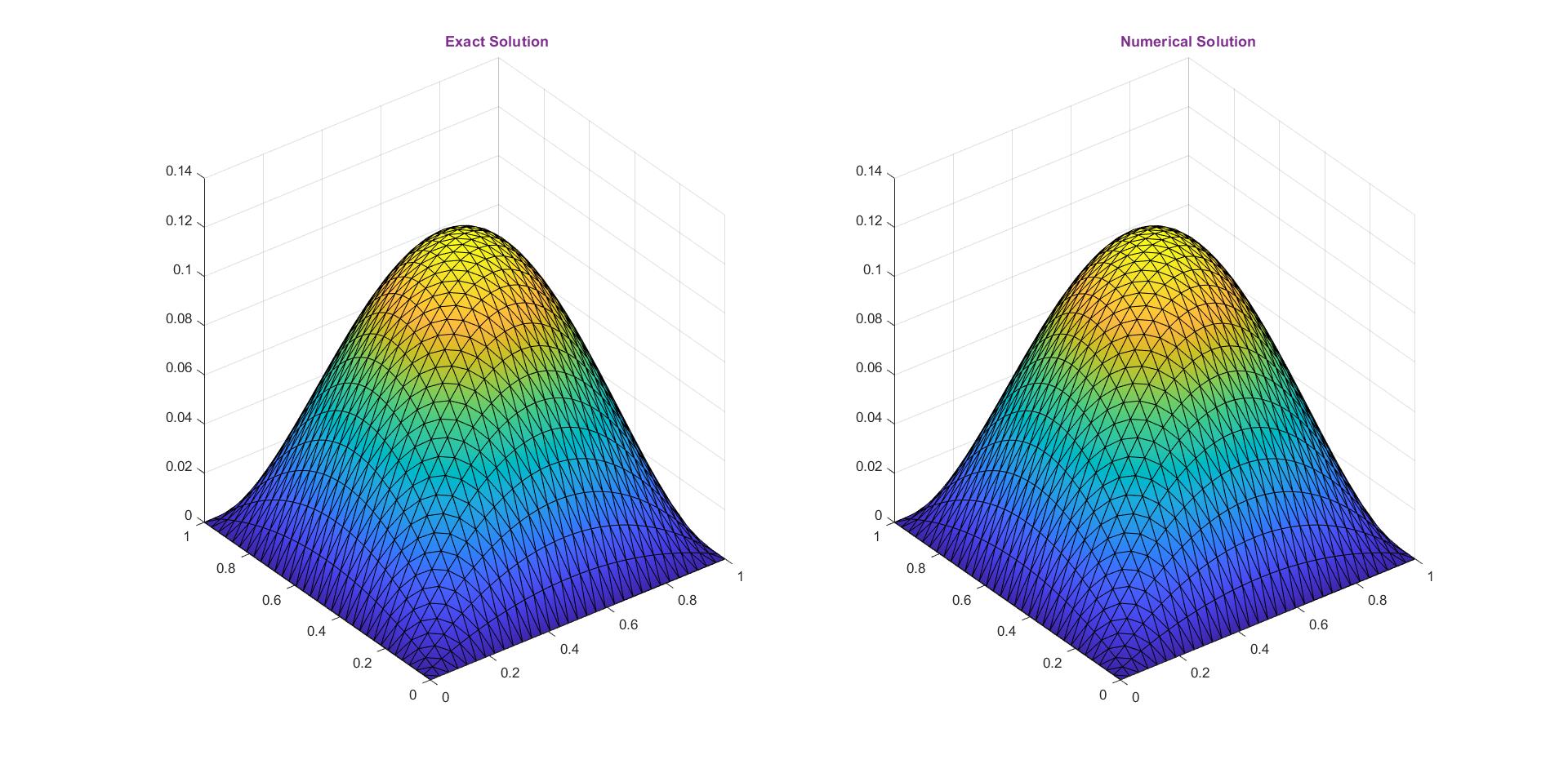, width=15 cm}
	\end{center}
    \caption{\emph{The exact and numerical solution at T=1,  $\alpha$ =0.5.}}
    \label{fig2}
\end{figure}

%%%%%%%%%%%%%%%%%%%%%%%%%%%%%%%%%%

%\section{Conclusions}
% to fill...

%%%%%%%%%%%%%%%%%%%%%%%%%%%%%%%%%%%
%%%%%%%%%%%%%%%%%%%%%%%%%%%%%%%%%%

\section{Acknowledgement}
 The authors acknowledge anonymous reviewers for many helpful suggestions and comments. Also, the authors would like to sincerely thank Professor Chaobao Huang for the valuable comments and suggestions.

%%%%%%%%%%%%%%%%%%%%%%%%%%%%%%%%%%%

%%%%%%%%%%%%%%%%%%%%%%%%%%%%%%%%%%

\begin{thebibliography}{20}
%\bibitem{[szy]} Z. Szyma$\acute{n}$ska, C. Morales-Rodrigo, M. Lachowicz, M. Chaplain, \emph{Mathematical modelling of cancer invasion of tissue: the role and effect of nonlocal interactions},  Math. Models Methods Appl. Sci., 19-2 (2009), 257-281.
 \bibitem{[me]} S. Chaudhary, \emph{Finite element analysis of nonlocal coupled parabolic problem using Newton’s method}, Comput. Math. Appl. 75-3 (2018), 981-1003.

\bibitem{[sk]} S. Chaudhary, V. Srivastava, V. V. K. Srinivas Kumar, B. Srinivasan, \emph{Finite element approximation of nonlocal parabolic problem}, Numer. Methods Partial Differ. Eq., 33 (2017) 786-313.

\bibitem{[ch]} M. Chipot, B. Lovat, \emph{On the asymptotic behaviour of some nonlocal problems}, Positivity 3 (1999), 65-81.

    \bibitem{[r6new]} S. B. Menezes, \emph{Remarks on weak solutions for a nonlocal parabolic problem}, Int. J. Math. Math. Sci., 2006 (2006), 1-10.

        \bibitem{[r4new]} S. Rida, A. El-Sayed, A. Arafa, \emph{Effect of bacterial memory dependent growth by using fractional derivatives reaction-diffusion chemotactic model}, J. Stat. Phys., 140-4 (2010), 797-811.

\bibitem{[r5new]} Y. Luchko, M. Rivero, J. Trujillo, M. Pilar Velasco, \emph{Fractional models, non-locality, and complex systems}, Comput. Math. Appl., 59-3 (2010), 1048-1056.

\bibitem{[rbl]} R. Robalo, R. Almeida, M. Coimbra, J. Ferreira, \emph{A reaction diffusion model for a class of nonlinear parabolic equations with moving boundaries: Existence, uniqueness, exponential decay and simulation}, Appl. Math. Model., 38-23 (2014), 5609-5622.

\bibitem{[hr12]} C. Huang, M. Stynes, \emph{Optimal spatial $H^1-norm$ analysis of a finite element method for a time-fractional diffusion equation},  J. Comput. Appl. Math., 367 (2020), 112435.

\bibitem{[ach1]} C. Huang, M. Stynes, \emph{Optimal  $H^1$ spatial convergence of a fully discrete finite element method for the time-fractional Allen-Cahn equation},  Adv. Comput. Math., 46-4 (2020), 63.

\bibitem{[kwn]} K. Mustapha, W. McLean, \emph{A second-order accurate numerical method for a fractional wave equation}, Numer. Math., 105 (2007), 481-510.

\bibitem{[r02]} M. Stynes, E. O'Riordan, J. Gracia, \emph{Error analysis of a finite difference method on graded meshes for a time-fractional diffusion equation}, SIAM J. Numer. Anal., 55 (2017), 1057-1079.

\bibitem{[r11]} N. Kopteva, \emph{Error analysis of the L1 method on graded and uniform meshes for a fractional derivative problem in two and three dimensions}, Math. Comp., 8 (2019) 2135-2155.

\bibitem{[mn1]} J. Manimaran, L. Shangerganesh, A. Debbouche and V. Antonov, \emph{Numerical solutions for time-fractional cancer invasion system with nonlocal diffusion}, Front. Phys., 7 (2019), 93.

\bibitem{[mn2]} J. Manimaran, L. Shangerganesh and A. Debbouche, \emph{Finite element error analysis of a time-fractional nonlocal diffusion equation with the Dirichlet energy},  J. Comput. Appl. Math., 382 (2021), 113066.

\bibitem{[mn3]} J. Manimaran, L. Shangerganesh, \emph{Error estimates for Galerkin finite element approximations of time-fractional nonlocal diffusion equation},  Int. J. Comput. Math., 98-7 (2020), 1365-1384.
	
\bibitem{[r1]} K. Diethelm, \emph{The Analysis of Fractional Differential Equations: An Application-Oriented Using Differential Operators of Caputo Type}, Lecture Notes in Mathematics, Springer, (2010).
	
%\bibitem{[r2]} H. Chen, Xiaohan Hu, J. Ren, T. Sun, Y. Tang, \emph{L1 Scheme on Graded Mesh for the Linearized Time Fractional KDV equation with Initial Singularity}, Int. J. Model. Simul. Sci. Comput., 10-1 (2018), 1941006.

\bibitem{[r3]} T. Gudi, \emph{Finite element method for a nonlocal problem of Kirchhoff type}, SIAM J. Numer. Anal., 50-2 (2012), 657-668.

\bibitem{[r5]}  H. Liao, W. Mclean, J. Zhang,  \emph{A Discrete Gr\"onwall Inequality with Applications to Numerical Schemes for Sub-diffusion Problems}, SIAM J. Numer. Anal., 57-1 (2019), 218-237.

\bibitem{[r5a]} J. Ren, H. Liao, J. Zhang, Z. Zhang,  \emph{Sharp H1-norm error estimates of two time-stepping schemes for reaction–subdiffusion problems}, J. Comput. Appl. Math., 389 (2021), 113352.

\bibitem{[r6]} M. Wheeler,  \emph{A Priori $l_2$ error estimates for Galerkin approximations to parabolic partial differential equations}, SIAM J. Numer. Anal., 10 (1973), 723-759.

\bibitem{[r7]} R. Rannacher and R. Scott, \emph{Some optimal error estimates for piecewise linear finite element approximation}, Math. Comp., 38 (1982), 490-507.
	
\bibitem{[vth]} V. Thom$\acute{e}$e, \emph{Galerkin Finite Element Methods for Parabolic Problems}, Second revised and expanded ed., Springer, Berlin, (2006).
	
\bibitem{[r08]}  K. Sakamoto, M. Yamamoto, \emph{Initial value/boundary value problems for fractional diffusion-wave equations and applications to some inverse problems}, J. Math. Anal. Appl., 382 (2011) 426-447.

\bibitem{[r8]} B. Jin, B. Li, Z. Zhou, \emph{Numerical analysis of nonlinear subdiffusion equations}, SIAM J. Numer. Anal., 56(1) (2018) 1-23.

     \bibitem{[xzw1]} X. Zheng, W.Hong, \emph{An Error Estimate of a Numerical  Approximation to a Hidden-Memory Variable-Order Space-Time Fractional Diffusion Equation}, SIAM J. Numer. Anal., 58-5 (2020), 2492-2514.
	
	   \bibitem{[xzw3]} X. Zheng, W.Hong, \emph{Optimal-order error estimates of finite element approximations to variable-order time-fractional diffusion equations without regularity assumptions of the true solutions}, IMA J. Numer. Anal., 41-2 (2020), 1522-1545.
	
	    \bibitem{[xzw4]} X. Zheng, W.Hong, \emph{A time-fractional diffusion equation with space-time dependent hidden-memory variable order: analysis and approximation}, Bit Numer. Math., 61 (2021), 1453-1481.
	
\bibitem{[r9]} M. Maskari, S. Karaa, \emph{Numerical approximation of semilinear subdiffusion equations with nonsmooth intial data}, SIAM J. Numer. Anal., 57 (2019) 1524-1544.

\bibitem{[r10]} D. Li, C. Wu, Z. Zhang, \emph{Linearized Galerkin FEMs for nonlinear time fractional parabolic problems with nonsmooth solutions in time direction}, J. Sci. Comput., 80 (2019) 403-419.

\bibitem{[r13]} D. Li, H. Qin, J. Zhang, \emph{Sharp pointwise-in-time error estimate of L1 scheme for nonlinear subdiffusion equations},  arXiv:2101.04554v1, (2021).

\bibitem{[r14]} S. Karaa, \emph{Galerkin type methods for semilinear time-fractional diffusion problems},  J. Sci. Comput. 83, 46 (2020).
	
\bibitem{[ddgm]} C. Huang, M. Stynes, \emph{A direct discontinuous Galerkin method for a time-fractional diffusion equation with a Robin boundary condition},  Appl. Numer. Math., 135 (2019), 15-29.

   %\bibitem{[lrs]}  M. G. Larson and F. Bengzon. \emph{The Finite Element Method: Theory,
%Implementation, and Applications}, Texts in Computational Science and
%Engineering, Springer, 2013.
\end{thebibliography}
\end{document}